\documentclass{amsart}
\title{Prikry On Extenders, Revisited}
\author{Carmi Merimovich}
\address{
The Academic College of Tel-Aviv
\\
4 Antokolsky St.
\\
Tel-Aviv 64044
\\
Israel
}
\email{carmi@mta.ac.il}
\date{May 8, 2003}


\usepackage{amssymb}

\usepackage{pb-diagram}

\theoremstyle{plain}
\newtheorem{theorem}{Theorem}[section]
\newtheorem*{theorem*}{Theorem}
\newtheorem*{Expected*}{Aim}
\newtheorem{lemma}[theorem]{Lemma}
\newtheorem{proposition}[theorem]{Proposition}
\newtheorem{claim}[theorem]{Claim}
\newtheorem{corollary}[theorem]{Corollary}

\theoremstyle{definition}
\newtheorem{definition}[theorem]{Definition}

\theoremstyle{remark}

\newcommand{\fA}{{\mathfrak{A}}}

\newcommand{\ga}{\alpha}
\newcommand{\gb}{\beta}
\newcommand{\gc}{\chi}

\newcommand{\gga}{\gamma}

\newcommand{\gk}{\kappa}
\newcommand{\gl}{\lambda}
\newcommand{\gm}{\mu}
\newcommand{\gn}{\nu}

\newcommand{\gp}{\pi}

\newcommand{\gr}{\rho}
\newcommand{\gs}{\sigma}
\newcommand{\gt}{\tau}

\newcommand{\gw}{\omega}
\newcommand{\gx}{\xi}

\newcommand{\gz}{\zeta}

\newcommand{\func}{\mathop{:}}
\newcommand{\VN}[1]{\Check{#1}}
\newcommand{\CN}[1]{\underset{\widetilde{}}{#1}}
\newcommand{\GN}[1]{\Dot{#1}}

\newcommand{\union}{\cup}
\newcommand{\bigunion}{\bigcup}
\newcommand{\intersect}{\cap}
\newcommand{\bigintersect}{\bigcap}
\newcommand{\forces}{\mathrel\Vdash}
\newcommand{\incompatible}{\perp}

\newcommand{\compatible}{\parallel}
\newcommand{\decides}{\mathrel\Vert}

\newcommand{\subelem}{\prec}

\newcommand{\append}{\mathop{{}^\frown}}
\newcommand{\restricted}{\mathrel{\restriction}}
\newcommand{\upto}{\mathord{<}}
\newcommand{\uptoeq}{\mathord{\leq}}

\newcommand{\gen}[1]{\mathfrak{g}(#1)}
\newcommand{\power}[1]{\lvert#1\rvert}

\newcommand{\On}{\ensuremath{\text{On}}}
\newcommand{\Cn}{\ensuremath{\text{Cn}}}
\newcommand{\ordered}[1]{\ensuremath{\langle #1 \rangle}}

\newcommand{\set}[1]{\ensuremath{\{ #1 \}}}
\newcommand{\setof}[2]{\ensuremath{\{ #1 \mid #2 \}}}
\newcommand{\ordof}[2]{\ensuremath{\ordered{ #1 \mid #2 }}}
\newcommand{\formula}[1]{{}^{\ulcorner} #1 {}^{\urcorner}}

\DeclareMathOperator{\geqE}{{\mathrel{\geq_{{E}}}}}
\DeclareMathOperator{\gtE}{{\mathrel{>_{{E}}}}}
\DeclareMathOperator{\ltE}{{\mathrel{<_{{E}}}}}

\DeclareMathOperator{\ltj}{{\mathrel{{<_{{j}}}}}}
\DeclareMathOperator{\gtj}{{\mathrel{{>_{{j}}}}}}

\DeclareMathOperator{\geqj}{{\mathrel{{\geq_{{j}}}}}}

\DeclareMathOperator{\ltT}{{\mathrel{{<_{{T}}}}}}

\DeclareMathOperator{\crit}{crit}
\DeclareMathOperator{\dom}{dom}
\DeclareMathOperator{\pcf}{pcf}

\DeclareMathOperator{\ot}{ot}
\DeclareMathOperator{\cf}{cf}
\DeclareMathOperator{\tcf}{tcf}
\DeclareMathOperator{\supp}{supp}
\DeclareMathOperator{\Ult}{Ult}
\DeclareMathOperator{\Suc}{Suc}
\DeclareMathOperator{\Lev}{Lev}

\DeclareMathOperator{\mc}{mc}

\newcommand{\Es}{{\ensuremath{\bar{E}}\/}}

\newcommand{\PE}{{\mathbb{P}_{E}\/}}
\newcommand{\PES}{{\mathbb{P}^*_{E}\/}}
\newcommand{\sZ}{\mathbb{Z}}

\begin{document}
\begin{abstract}
We present a modification to the Prikry on Extenders forcing
notion allowing the blow up of the power set of a large cardinal,
change its cofinality to $\gw$ without adding bounded subsets,
 working directly from  \emph{arbitrary} extender
(e.g., $n$-huge extender).

Using this forcing, starting from a superstrong cardinal $\gk$,
we construct a model in which the added Prikry sequences are
a scale in the normal Prikry sequence.
\end{abstract}
\maketitle
\section{Introduction}
In \cite{PrikryExtender}, Gitik and Magidor introduced technology
to blow up the power of a large cardinal, $\gk$, 
change its cofinality
to $\gw$, while preserving all cardinals and adding no new bounded
subsets to $\gk$. As defined this forcing notion could utilize extenders
of length $\gl$ derived from an elementary embedding $j \func V \to M$ satisfying
$\crit j = \gk$, $M \supseteq M^\gk$ as long as 
                $\gl \leq \sup_{f \func \gk \to \gk} j(f)(\gk)$.
If the elementary embedding is strong enough and 
$\gl > \sup_{f\func \gk \to \gk} j(f)(\gk)$ is demanded, a preparation forcing
is done adding a generic function $f$ satisfying $j(f)(\gk) > \gl$
(of course, after lifting $j$ to the generic extension).

In this work we present a modification
of the Gitik-Magidor forcing
which  allows to work directly with whatever
extender is presented. For example, this modification allows us to
start from $\gk$ which is $\gl$-supercompact, (i.e., there is
$j \func V \to M \supset M^\gl$, $\crit(j) = \gk$) and get
a generic extension in which $2^{\gk} = \power{j(\gl)}$ and all originally
regular cardinals in the range $[\gk,j(\gl)]$ change their cofinality
to $\gw$. In fact substituting $j_n(\gk)$ for $\gl$ in the above supercompact
embedding (that is starting from an $n$-huge cardinal)
we can get  a generic extension in which $2^\gk = \power{j_{n+1}(\gk)}$,
$\forall \gm \in [\gk, j_n(\gk)]_{\text{Reg}}$ $\cf \gm = \gw$.

As an added benefit we get more closure of the Prikry order. Namely,
in the original forcing the Prikry order is $\gk$-closed while
in the modified forcing it is $\gk^+$-closed.

The structure of this work is as follows: In section \ref{ExtenderIndexing}
we define what is an extender in this work. In section \ref{PEForcing}
we give a detailed presentation of the modified Gitik-Magidor
forcing notion 
assuming the extender we use is at most superstrong.
With the exception of the Prikry property proof,
up to \ref{NewStuff} 
 we give adaptation
of the proofs from \cite{PrikryExtender} to the modified forcing.
This section  culminates with
\begin{theorem*}[\ref{PrikryExtenderMain}]
Assume GCH, $j \func V \to M \supset M^\gk$, $\crit(j)=\gk$.
Then there is a cardinal
preserving generic extension
in which $2^\gk = \power{j(\gk)}$, $\cf \gk = \gw$, and there are no new bounded
subsets of $\gk$.
\end{theorem*}
Section \ref{pcfApplication} is a result of a talk with Gitik where
 he pointed out that if there is no restriction on the extender then
we can start from a superstrong  cardinal and get a sequence of functions
which are unbounded in the product of the normal Prikry sequence.
That is
\begin{theorem*}[\ref{ScaleTheorem}]
Assume GCH,
        $j\func V \to M \supset M^\gk$, $M \supset V_{j(\gk)}$,
 $\crit{j} = \gk$.
Then there is a cardinal preserving generic extension in which 
$\cf \gk = \gw$,
$\gk^\gw = j(\gk)$, and
$\forall \gl \in [\gk, j(\gk))_{\text{Reg}}$ there is a function
$G^{\gl}\func \gw \to \gk$ such that
$\tcf \prod G^{\gl}/D = \gl$,
and
$\ordof {G^{\gl}/D} {\gl \in [\gk,j_E(\gk))_{\text{Reg}}}$ is a scale in
        $\prod G^{\gk}$,
where $D$ is the cofinite filter over $\gw$.
\end{theorem*}
In section \ref{GenericByIteration} we show
several ways (mainly because we do not know 
the `right' way, if it exists at all)
to generate a generic filter in $V$ over
the $\gw$ iterate of $V$.
In section \ref{PEForcingGeneral} we define the forcing notion
for stronger elementary embeddings allowing us to get:
\begin{theorem*}[\ref{PrikryExtenderMainGeneral}]
Assume GCH, $j \func V \to M \supset M^\gl$, $\crit (j)=\gk$.
Then there is 
a generic extension
in which $2^\gk = \power{j(\gl)}$, 
	all cardinal up to $\gk$ and above $\gl$ are preserved,
	$\forall \gm\in [\gk,\gl]_{\text{Reg}}$
		$\cf \gm = \gw$, and there are no new bounded
subsets of $\gk$.
\end{theorem*}
We do not give proofs in this section since they are essentially the same
proofs as in section \ref{PEForcing}.

This work is largely self contained. Of course knowledge of \cite{PrikryExtender}
will make the reading very easy.
The notation we use is standard. We assume fluency with 
forcing (as in say, \cite{KunenBook}),
large cardinals and extenders (as in \cite{AkiBook}),
and some basic $\pcf$ theory (as in \cite{HolzBook}).
\section{Elementary embeddings and Extenders.} \label{ExtenderIndexing}
\begin{definition}
Let $j\func V \to M$ be an elementary embedding, $\crit(j)=\gk$.
\begin{enumerate}
\item
	The generators\footnote{The definition of generators differs slightly
	from the usual one since we use only one  ordinal to index our extenders.}
	 of $j$ are defined by induction as
	\begin{align*}
	& \gk_0 = \crit(j),
	\\
	& \gk_\gx = \min \setof {\gl \in \On} {\forall \gx'<\gx\ 
			\forall \gm \in \On\ 
			\forall f\func \gm \to \On\ 
			j(f)(\gk_{\gx'}) \neq \gl}.
	\end{align*}
	If the induction terminates, then we have a set of generators for $j$:
	\begin{align*}
	& \gen{j} = \setof {\gk_\gx} {\gx < \gx^*}.
	\end{align*}
\item
	For $\ga, \gb \in \On$ we say $\ga \ltj \gb$ if
	\begin{enumerate}
	\item
	        $\ga < \gb$.
	\item
	        There is $f\func \gm \to \On$ such that
	                $j(f)(\gb) = \ga$.
	\end{enumerate}
\item
	Assume $\ga \in \On$, $\gl \in \Cn$ is minimal such that
		$j(\gl) > \ga$. We set
	\begin{align*}
	\forall A \subseteq \gl\ A \in E(\ga) \iff \ga \in j(A).
	\end{align*}
	It is well know that $E(\ga)$ is a $\gk$-complete ultrafilter
	over $\gl$. Note that for a cardinal $\gm$,
	an ultrafilter generated by 
	$j(\gm) < \ga < \bigunion j''\gm^+$ is isomorphic to an ultrafilter
	on $\gm$ generated by some $\gb < j(\gm)$.
	
\end{enumerate}
\end{definition}
\begin{definition}
Let $j\func V \to M \supset M^\gl$ be an elementary embedding, $\crit(j)=\gk$,
	$\gen{j} \subset j(\gl)$. The extender $E$ derived from $j$
	is the system
	\begin{align*}
	E = \ordered{\ordof{E(\ga)} {\ga \in j(\gl) \setminus j''\gl},
			 \ordof{\gp_{\gb,\ga}} {\ga,\gb \in j(\gl)\setminus j''\gl,\ 
						\ga \ltj \gb}
			}.
	\end{align*}
	where
	\begin{enumerate}
	\item
		$\forall A \subseteq \gk\ A \in E(\ga) \iff \ga \in j(A)$.
	\item
		For $\ga,\gb \in j(\gl) \setminus j''\gl$, $\ga \ltj \gb$,
		 the function $\gp_{\gb,\ga}\func \gl \to \gl$ is such that
		$j(\gp_{\gb,\ga})(\gb) = \ga$.
		(Note that $\ga \ltj \gb$ means there are many such functions.
		Any one of them will do as $\gp_{\gb,\ga}$.)
	\end{enumerate}
	We assume it is known that $j$ can be reconstructed from $E$,
	i.e., $j$ is the canonical embedding $j \func V \to \Ult(V,E)$.
	We will use $\ltE$, $\dom E$ as synonyms for $\ltj$,
			$j(\gl) \setminus j''\gl$, respectively.
\end{definition}
\begin{claim}
Assume $\gl^{\upto \gl} = \gl$, $j\func V \to M$, $\crit(j) = \gk$, $M \supseteq M^\gl$. Then
$\ltj \restricted j(\gl)$ is $\gl^+$-directed.
\end{claim}
\begin{proof}
Let $X \in [j(\gl)]^{\leq \gl}$. We need to find $\gb < j(\gl)$ such that
$\forall \ga \in X$ $\gb \gtj \ga$.

Let us fix a function $e\func \gl \xrightarrow{\text{onto}} [\gl]^{\upto \gl}$ such that
for each $A \in [\gl]^{\upto \gl}$, $e^{-1}A$
is unbounded in $\gl$.
Of course, $j(e)\func j(\gl) \xrightarrow{\text{onto}} 
		([j(\gl)]^{\upto j(\gl)})_M$. 

Let $\gm = \sup X$. Since $X \in M$ we get $\gm < j(\gl)$.
Since $X \in ([j(\gl)]^{\upto j(\gl)})_M$ there are $\gb \geq \gm$, $g$ such that
$j(e)(\gb) = X$, $j(g)(\gb)=j''\ot(X)$. 
We show that $\gb \geqj \ga$ for all $\ga \in X$. So,
let $\ga \in X$.

We let $\gx = \ot(X \intersect \ga)$.
Then we set $\forall \gn < \gl$ $g_\gx(\gn) = \ot(g(\gn) \intersect \gx)$.
Thus $j(g_\gx)(\gb) = \ot(j(g)(\gb) \intersect j(\gx)) =
	\ot(j''\ot(X) \intersect j(\gx)) = \gx$.
We set $\forall \gn < \gl$ $f(\gn) = min \setof {\gga \in e(\gn)} 
	{\ot(e(\gn) \intersect \gga) = g_\gx(\gn))}$. Then
$j(f)(\gb) = \min \setof {\gga \in X} {\ot(X \intersect \gga) = \gx}= \ga$.
\end{proof}
We give  basic definitions for iterating elementary embedding, and
a proposition regarding their representation.
\begin{definition}
Assume $j\func V \to M$ is an elementary embedding. We define by induction
$\forall n < \gw$
\begin{align*}
	& j_{0,1} = j,\ M_0 = V,
	\\
	& j_{n+1,n+2} = j(j_{n,n+1}) \func M_{n+1} \to M_{n+2}.
\end{align*}
We `complete' the list of $j$'s by setting $\forall n < m < \gw$
	\begin{align*}
	& j_{n,m} = j_{m-1,m} \circ \dotsb \circ j_{n, n+1},
	\\
	& j_n = j_{0,n}.
	\end{align*}
\end{definition}
\begin{proposition}
Assume $j\func V \to M \supset M^\gl$, 
$\crit(j) = \gk$,
$\gen{j} \subset j(\gl)$.
$\gt \in j_n(j(\gl) \setminus j''(\gl))$. Then there is 
        $\gt^* \in j(\gl)$ such that $j_n(\gt^*) \gtj \gt$.
\end{proposition}
\begin{proof}
The proof is by induction on $n$.
\begin{itemize}
\item
        $n=1$:
        We choose $f \func \gl \to \gl$, $\ga \in j(\gl) \setminus j''\gl$
	 such that
        $j(f)(\ga) = \gt$. Since $\ltj \restricted j(\gl)$ is 
	$\gl^+$-directed, there is
        $\gt^* \in j(\gl)$ such that 
		$\forall \gn < \gl$ $\gt^* \gtj f(\gn)$.
        Hence $j(\gt^*) \gtj j(f)(\ga) = \gt$.
\item
        $n > 1$:
        By induction there is $\gt^{\prime*} \in j_{n-1}(j(\gl))$
        such that $j_{n-1,n}(\gt^{\prime*}) \gtj \gt$. By induction, 
        there is
        $\gt^{*} \in j(\gl)$ such that 
	$j_{n-1}(\gt^{*}) \gtj \gt^{\prime*}$.
        So
        \begin{align*}
        j_n(\gt^{*}) = j_{n-1,n}(j_{n-1}(\gt^{*})) \gtj
                j_{n-1,n}(\gt^{\prime*}) \gtE \gt.
        \end{align*}
\end{itemize}
\end{proof}
\begin{corollary}
Let $n < \gw$ and $M_n$ be the $n$-th $j$-iterate of $V$ and $x \in M_n$. Then
there are $f\func [\gl]^n \to V$ and $\ga \in j(\gl)$ such that
	$x = j_n(f)(\ga, j(\ga), \dotsc, j_{n-1}(\ga))$.
\end{corollary}
\section{$\PE$-Forcing} \label{PEForcing}
In this section we give a detailed presentation of 
the Prikry on Extenders
forcing notion, due to Gitik and Magidor \cite{PrikryExtender},
 assuming $\gk$ carries a super-strong extender at most.
Hence we begin with definitions of trees and functions on trees which
are essentially on $\gk$.
\begin{definition}
Assume $T \subseteq [\gk]^{\upto \gw}$. 
For
        $\ordered{\gm_0, \dotsc, \gm_k},
                \ordered{\gn_0, \dotsc, \gn_n} \in T$
we define
        $\ordered{\gm_0, \dotsc, \gm_k} \ltT 
                \ordered{\gn_0, \dotsc, \gn_n}$
if
\begin{enumerate}
\item
        $k < n$.
\item
        $\ordered{\gm_0, \dotsc, \gm_k} =
                \ordered{\gn_0, \dotsc, \gn_k}$.
\end{enumerate}
Clearly $\ordered{T, \ltT}$ is a tree.
We always assume that $\ordered{} \in T$.
\end{definition}
\begin{definition}
Assume $T \subseteq [\gk]^{\upto \gw}$ is a tree. Then
\begin{enumerate}
\item
        $\Suc_T(\ordered{\gn_0, \dotsc, \gn_k}) =
                 \setof {\gn < \gk}
                {\ordered{\gn_0, \dotsc, \gn_k, \gn} \in T}$.
\item
        $\forall k<\gw$  $\Lev_k(T) = T \intersect [\gk]^k$.
\end{enumerate}
Note that $\ordered{} \in T$ implies $\:\Lev_0(T) = \set{\ordered{}}$.
\end{definition}
\begin{definition}
Assume $T \subseteq [\gk]^{\upto \gw}$ is a tree, 
        $\ordered{\gn_0,\dotsc, \gn_k} \in T$. Then
\begin{align*}
T_{\ordered{\gn_0, \dotsc, \gn_k}} = 
        \setof {\ordered{\gn_{k+1}, \dotsc, \gn_n} \in [\gk]^{\upto \gw} }
                {\ordered{\gn_0, \dotsc, \gn_k, \gn_{k+1}, \dotsc, \gn_n} \in T}.
\end{align*}
\end{definition}
\begin{definition}
Assume $T \subseteq [\gk]^{\upto \gw}$ is a tree, $A \in [\gk]^k$. Then
\begin{align*}
T\restricted A =
        \setof {\ordered{\gn_{1}, \dotsc, \gn_n} \in T}
                {n < \gw,\ \ordered{\gn_0, \dotsc, \gn_k} \in A}.
\end{align*}
\end{definition}
\begin{definition}
Assume $T^\gx \subseteq [\gk]^{\upto \gw}$ is a tree for all $\gx < \gl$.
Then $T = \bigintersect_{\gx < \gl} T^\gx$ is defined by induction on $k$ as:
\begin{enumerate}
\item
        $\Lev_0(T) = \set{\ordered{}}$.
\item
        $\ordered{\gn_0, \dotsc, \gn_k} \in T \implies$
                $\Suc_T({\ordered{\gn_0, \dotsc, \gn_k}}) =
                        \bigintersect_{\gx < \gl} 
                               \Suc_{T^\gx}(\ordered{\gn_0, \dotsc, \gn_k})$.
\end{enumerate}
\end{definition}
\begin{definition}
Let $F$ be a function such that $\dom F \subseteq [\gk]^{\upto \gw}$
is a tree and
        $\ordered{\gn_0,\dotsc, \gn_k} \in \dom F$. Then
        $F_{\ordered{\gn_0, \dotsc, \gn_k}}$ is a function such that
\begin{enumerate}
\item
        $\dom (F_{\ordered{\gn_0, \dotsc, \gn_k}}) = 
                (\dom F)_{\ordered{\gn_0, \dotsc, \gn_k}}$.
\item
        $F_{\ordered{\gn_0, \dotsc, \gn_k}} (\gn_{k+1}, \dotsc, \gn_n) = 
                F(\gn_0, \dotsc, \gn_k, \gn_{k+1}, \dotsc, \gn_n)$.
\end{enumerate}
\end{definition}
From now on we assume GCH and the existence of
 $j \func V \to M \supset M^\gk$,
	$\crit(j) = \gk$, $\gen{j} \subset j(\gk)$, and $E$ is the
extender derived from $j$. Recall $\dom E = j(\gk) \setminus \gk$.
\begin{definition}
$T \subseteq [\gk]^{\upto \gw}$ is $E(\ga)$-tree if
        $\forall \ordered{\gn_0, \dotsc, \gn_k} \in T$
                $\Suc_T(\ordered{\gn_0, \dotsc, \gn_k})
                                \in E(\ga)$.
\end{definition}
We recall the definition of filter product in order to define powers
of $E(\ga)$.
\begin{definition}
We define powers of $E(\ga)$ by induction as follows:
\begin{enumerate}
\item
        For $k = 0$: $E^0(\ga) = \set{\emptyset}$.
\item
        For $k > 0$:
                $\forall A \subseteq [\gk]^k$ $A \in E^{k}(\ga) \iff$
        \begin{multline*}
        \setof{\ordered{\gn_0,\dotsc, \gn_{k-1}} \in [\gk]^{k-1}} 
        {
\\
        \setof {\gn_{k} \in \gk} 
                {\ordered{\gn_0, \dotsc, \gn_{k-1}, \gn_k} \in A} \in E(\ga)
        } \in E^{k-1}(\ga).
\end{multline*}
\end{enumerate}
Note that $E^1(\ga) = E(\ga)$. Recall that $E^k(\ga)$ is $\gk$-closed
ultrafilter on $[\gk]^k$.
\end{definition}
The following is straightforward.
\begin{proposition}
Assume $T \subseteq [\gk]^{\upto \gw},\  T^\gx \subseteq [\gk]^{\upto \gw}$ are
         $E(\ga)$-trees for all $\gx < \gl$, $\gl < \gk$. Then
\begin{enumerate}
\item
        $\forall k<\gw\ \Lev_k(T) \in E^k(\ga)$.
\item
        $\ordered{\gn_0, \dotsc, \gn_k} \in T \implies$
        $T_{\ordered{\gn_0, \dotsc, \gn_k}}$ is $E(\ga)$-tree.
\item
        $A \in E^k(\ga) \implies$
                $T \restricted A$ is $E(\ga)$-tree.
\item
        $\bigintersect_{\gx < \gl} T^{\gx}$ is $E(\ga)$-tree.
\end{enumerate}
\end{proposition}
We are ready to present the forcing notion:
\begin{definition} \label{PEdefinitionK}
A condition $p$ in $\PE$ is of the form
        $\ordered{f, \ga, F}$
where
\begin{enumerate}
\item
    $f\func d \to [\gk]^{\upto \gw}$ is such that
\begin{enumerate}
        \item
                $d \in [\dom E]^{\uptoeq \gk}$.
        \item
                $\forall \gb \in d$  $\ga \geqE \gb$.
\end{enumerate}
\item
        $F \func T \to [d]^{\upto \gk}$ is such that
\begin{enumerate}
        \item
                $T$ is $E(\ga)$-tree.
         \item  \label{ForceReqSupport}
                $\forall \ordered{\gn_0, \dotsc, \gn_k, \gn} \in T$
                        $F_{\ordered{\gn_0, \dotsc, \gn_k}}(\gn) \subseteq
                                        d$. 
        \item
                $\forall \ordered{\gn_0, \dotsc, \gn_k} \in T$
                        $j(F_{\ordered{\gn_0, \dotsc, \gn_k}})(\ga) = 
                                                                   j''d$.
        \item   \label{ForceReqMinimalAlways}
                $\forall \ordered{\gn_0, \dotsc, \gn_k, \gn} \in T$
                        $\gk \in F_{\ordered{\gn_0, \dotsc, \gn_k}}(\gn)$.
                Note: This condition implies $\gk \in d$!
        \item  \label{ForceReqMonotone}
                $\forall \gb \in d$ $\forall \ordered{\gn_0, \dotsc, \gn_k} 
                                                        \in T\ $
                \begin{align*}
                f(\gb) \append \ordof {\gp_{\ga,\gb}(\gn_i)} {i \leq k,\,
                        \gb \in F(\gn_0, \dotsc, \gn_i)}
                        \in [\gk]^{\upto \gw}.
                \end{align*}
\end{enumerate}
\end{enumerate}
We write $\supp p$, $\mc(p)$, $f^p$, $F^p$, $\dom p$ for
         $d$,       $\ga$,    $f$,   $F$, $T$.
\end{definition}
%\fi
%
%
%
\begin{definition} \label{PrikryBasic1}
Let $p, q \in \PE$. We say that $p$ is a Prikry extension of $q$
($p \leq^* q$ or $p \leq^0 q$) if
\begin{enumerate}
\item
        $\supp p \supseteq \supp q$.
\item
        $f^p \restricted \supp q = f^q$.
\item 
        $\dom p \subseteq \gp_{\mc(p), \mc(q)}^{-1} \dom q$.
\item
        \label{PrikryComute1}
        $\forall \ordered{\gn_0,\dotsc,\gn_k, \gn} \in \dom p$
        $\forall \gb \in F^q \circ \gp_{\mc(p), \mc(q)}
                        {}_{\ordered{\gn_0, \dotsc, \gn_k}}(\gn)$
        \begin{align*}
		\gp_{\mc(p), \gb}(\gn) =
			\gp_{\mc(q), \gb}(
				\gp_{\mc(p), \mc(q)}(\gn)).
        \end{align*}
\item
        $\forall \ordered{\gn_0,\dotsc,\gn_k} \in \dom p$
                $F^p(\gn_0, \dotsc, \gn_k) \supseteq
                        F^q \circ \gp_{\mc(p),\mc(q)}
                         (\gn_0, \dotsc, \gn_k)$.
\item
        $\forall \ordered{\gn_0,\dotsc,\gn_k} \in \dom p$
                $F^p(\gn_0, \dotsc, \gn_k) \setminus
                        F^q \circ \gp_{\mc(p),\mc(q)}
                         (\gn_0, \dotsc, \gn_k) \subseteq 
                                \supp p \setminus \supp q$.
\end{enumerate}
\end{definition}
Note that we can do without the requirements
\ref{ForceReqSupport}, \ref{ForceReqMinimalAlways}, \ref{ForceReqMonotone}
in definition \ref{PEdefinitionK}.
That is, if we define a forcing notion $P'_{\Es}$ in the same way
we defined $\PE$
but without these requirements then, using the above definition
and \ref{ArbitraryInSupport}, $\PE$ is $\leq^*$-densely
embeddable into $P'_{\Es}$.
\begin{definition} \label{dfn:enlarge1}
Let $q \in \PE$ and $\ordered{\gn} \in \dom q$. 
We define 
$q_{\ordered{\gn}} \in \PE$ to be $p$ where
\begin{enumerate}
\item
        $\supp p = \supp q$.
\item
        $\forall \gb \in \supp p$
	$f^p(\gb)=$
        $\begin{cases}
        f^q(\gb) \append \ordered{\gp_{\mc(q), \gb}(\gn)}   &
                \gb \in F^q(\gn).
        \\
	f^q(\gb)   &
                \gb \notin F^q(\gn).
        \end{cases}$
\item
        $\mc(p) = \mc(q)$.
\item
        $F^p = F^q_{\ordered{\gn}}$.
\end{enumerate}
\end{definition}
When we write $q_{\ordered{\gn_0, \dotsc, \gn_k}}$ we mean
        $(\dotsb (q_{\ordered{\gn_0}})_{\ordered{\gn_1}}\dotsb)_
                        {\ordered{\gn_k}}$.
\begin{definition}
Let $p,q \in \PE$. We say that $p$ is a $1$-point extension of $q$
($p \leq^1 q$) if
there is $\ordered{\gn} \in \dom q$ such that
        $p \leq^* q_{\ordered{\gn}}$.
\end{definition}
%
%
%
% n-point extension
%
\begin{definition}
Let $p,q \in \PE$. We say that $p$ is an $n$-point extension of $q$
($p \leq^n q$) if there are $p^n, \dotsc, p^0$ such that
\begin{align*}
        p=p^n \leq^1 \dotsb \leq^1 p^0=q.
\end{align*}
\end{definition}
\begin{definition}
Let $p,q \in \PE$. We say that $p$ is an extension of $q$
($p \leq q$) if there is $n$ such that $p \leq^n q$.
\end{definition}
While we have made a notational change from the original \cite{PrikryExtender} definition,
we hope it is obvious it is almost exactly the same forcing. The \emph{real}
changes from \cite{PrikryExtender} are:
\begin{enumerate}
\item
        We eliminated from the forcing definition the requirement
        that for each $p \in \PE$
        \begin{align*}
        \forall \gn \in \dom p\ 
        \power{\setof{\gb \in \supp p} {f^p(\gb) < \gp_{\mc(p), \gk}(\gn) }} \leq 
                                        \gp_{\mc(p), \gk} (\gn).
        \end{align*}
\item
        We  changed the rule for extending $f^q(\gb)$ in the definition
        of $ p = q_{\ordered{\gn}}$. It was
        \begin{align*}
                f^p(\gb) = 
                \begin{cases}
                        f^q(\gb) \append \gp_{\mc(q), \gb}(\gn) &
                                \gp_{\mc(q), \gk}(\gn) > \max f^q(\gb).
                \\
                        f^q(\gb)                &
                                \gp_{\mc(q), \gk}(\gn) \leq \max f^q(\gb).
                \end{cases}
        \end{align*}
        Now we have  the function $F$ in order to decide where to add
        a point:
        \begin{align*}
                f^p(\gb) = 
                \begin{cases}
                        f^q(\gb) \append \gp_{\mc(q), \gb}(\gn) &
                                \gb \in F^q(\gn).
                \\
                        f^q(\gb)                &
                                \gb \notin F^q(\gn).
                \end{cases}
        \end{align*}
\end{enumerate}
Putting it simply, originally the Prikry sequences had some constraint
because they carried also the information to which coordinate to add
points. We lift the limitation on the Prikry sequences by putting
this information in a separate place.

The following is immediate from the definition of the forcing notion.
\begin{claim} \label{CommuteOrder}
Let $p, q \in \PE$ such that $p \leq q$. Then
\begin{enumerate}
\item
        There is $\ordered{\gn_0, \dotsc, \gn_k} \in \dom q$ 
        such that
                $p \leq^* q_{\ordered{\gn_0, \dotsc, \gn_k}}$.
\item
        There are $r \in \PE$, $\ordered{\gn_0, \dotsc, \gn_k} \in \dom r$ 
        such that  $r \leq^* q$,
                $p = r_{\ordered{\gn_0, \dotsc, \gn_k}}$.
\end{enumerate}
\end{claim}
\begin{corollary}
Let $p \in \PE$. Then
\begin{enumerate}
\item
        Let $k < \gw$. Then
        $\setof {p_{\ordered{\gn_0, \dotsc, \gn_k}}}
                 {\ordered{\gn_0, \dotsc, \gn_k}
                             \in \dom p}$ is a maximal
                                anti-chain below $p$.
\item
        Let $k < \gw$, $\gs$ a formula in the forcing language.
        Assume $\forall\ordered{\gn_0, \dotsc, \gn_k} \in \dom p$
                $p_{\ordered{\gn_0, \dotsc, \gn_k}} \forces \gs$.
        Then $p \forces \gs$.
        
\end{enumerate}
\end{corollary}
\begin{definition}
Let $G$ be $\PE$-generic. Then
\begin{align*}
\forall \ga \in \dom E\ 
        G^\ga = \bigunion \setof {f^p(\ga)} {p \in G,\ \ga \in \supp p}.
\end{align*}
We write $\CN{G}^\ga$ for the $\PE$-name of $G^\ga$.
\end{definition}
\begin{lemma} \label{AnyMC}
\label{CommutShrink}
Let $q \in \PE$, $\ga \gtE \mc(q)$. Then there is $p \leq^* q$ with
$\mc(p) = \ga$.
\end{lemma}
\begin{proof}
We observe that for each $\gb \in \supp q$
\begin{align*}
	& j(\gp_{\ga, \gb})(\ga) = \gb,
\\
	& j(\gp_{\mc(q), \gb})(
				j(\gp_{\ga, \mc(q)})(\ga))  = \gb.
\end{align*}
So in $M$ we have
\begin{align*}
\forall \gb \in j'' \supp q \ 
		j(\gp)_{j(\ga), \gb}(\ga) =
			j(\gp)_{j(\mc(q)), \gb}(
				j(\gp_{\ga, \mc(q)})(\ga)).
\end{align*}
Let us set $F' = F^q \circ \gp_{\ga, \mc(q)}$. Then
\begin{align*}
\forall \ordered{\gn_0, \dotsc, \gn_k} \in \dom F'\ 
        j(F'_{\ordered{\gn_0, \dotsc, \gn_k}})(\ga) = j'' \supp q.
\end{align*}
We see that 
\begin{multline*}
\forall \ordered{\gn_0, \dotsc, \gn_k} \in \dom F'\ 
        A_{\ordered{\gn_0, \dotsc, \gn_k}} = \setof {\gn < \gk}
             {\forall \gb \in F'_{\ordered{\gn_0, \dotsc, \gn_k}}
                                                (\gn)\ 
\\
                \gp_{\ga,\gb}(\gn) =
                 \gp_{\mc(q), \gb}(\gp_{\ga, \mc(q)}(\gn)) }
                \in
                E(\ga).
\end{multline*}
So, let $F$ be $F'$ shrunken to these sets, namely
\begin{multline*}
\forall \ordered{\gn_0, \dotsc, \gn_k} \in \dom F\ 
                \Suc_{\dom F}  (\ordered{\gn_0, \dotsc, \gn_k}) = 
        \\
                \Suc_{\dom F'} (\ordered{\gn_0, \dotsc, \gn_k}) \intersect
                        A_{\ordered{\gn_0, \dotsc, \gn_k}}.
\end{multline*}
We set $p = \ordered{f^q, \ga, F}$ and by the above shrinkage we get 
$p \leq^* q$.
\end{proof}
Later on we do pull up of the form $\gp_{\ga, \gb}^{-1}$ freely.
It is implicitly assumed the shrinkage done in the above
proof, from $F'$ to $F$, is in this pull up.
\begin{proposition} \label{ArbitraryInSupport}
Let $q \in \PE$, $\ga \in \dom E$. Then there is $p \leq^* q$ with
$\ga \in \supp p$.
\end{proposition}
\begin{proof}
If $\ga \in \supp q$ then there is nothing to do, we set $p = q$.

If $\ga \notin \supp q$ we choose $\gga \gtE \ga$ such that 
        $\gga \gtE \mc(q)$.
By \ref{AnyMC}, there is $p' \leq^* p$ with $\mc(p') = \gga$.
We set $F$ as 
\begin{align*}
\forall \ordered{\gn_0, \dotsc, \gn_k, \gn} \in \dom F^{p'}\ 
        F_{\ordered{\gn_0, \dotsc, \gn_k}}(\gn) = 
        F^{p'}_{\ordered{\gn_0, \dotsc, \gn_k}}(\gn) \union \set {\ga},
\end{align*}
and then we set $p^* = \ordered{f^{p'} \union 
        \set{\ordered{\ga, \emptyset}}, \gga, F}$.
\end{proof}
We did not have to add $\ga$ to all values of $F$. Adding $\ga$
on an $E(\gga)$-large set would have been enough. This fact will be
used later on.

From the above propositions we see that 
for all $\ga \in \dom E$, $G^\ga$ is not empty.
In fact using density arguments we get:
\begin{proposition}
Let $G$ be $\PE$-generic. Then in $V[G]$:
\begin{enumerate}
\item $\ot G^\ga = \gw$.
\item $G^\ga$ is unbounded in $\gk$.
\item $\ga \not= \gb \implies G^\ga \not= G^\gb$.
\end{enumerate}
\end{proposition}
\begin{proposition}
$\PE$ satisfies $\gk^{++}$-cc.
\end{proposition}
\begin{proof}
Let $X \subseteq \PE$, $\power{X} = \gk^{++}$.
Since for each $p \in X$ we have $\power{\supp p} \leq \gk$, we can assume
that $\setof {\supp p} {p \in X}$ forms a $\Delta$-system. That is,
there is $d$ such that $\forall p,q \in X$ $\supp p \intersect \supp q = d$.
Since $\power{d} \leq \gk$ we have 
$\power{\setof{f} {f\func d \to [\gk]^{\upto\gw}}} \leq \gk^+$, 
so we can assume that
$\forall p,q \in X$ $\forall \gb \in d$ $f^p(\gb) = f^q(\gb)$.

Let us fix $2$ conditions, $p, q \in X$. Let $f = f^p \union f^q$
Then $f\func \supp p \union \supp q \to [\gk]^{\upto \gw}$.
Choose $\ga \gtE \mc(p), \mc(q)$.
Set $F' = F^p \circ \gp_{\ga, \mc(p)} \union F^q \circ \gp_{\ga, \mc(q)}$.
Let  $r' = \ordered{f, \ga, F'}$. Then $r'$ almost $\leq p,q$ in the
sense of \ref{CommutShrink}. Hence shrinking $\dom r'$
using \ref{CommutShrink} once with respect to $p$ and once with respect to
$q$  yields $r \leq p,q$.
\end{proof}
Up to this point we know that in a $\PE$-generic extension we have
\begin{enumerate}
\item
        $\cf \gk = \gw$.
\item
        $2^\gk = \power{j(\gk)}$.
\item
        Cardinals above $\gk^+$ are not collapsed.
\end{enumerate}
In order to see that no damage happens below $\gk$ we use the Prikry
ordering.
\begin{proposition} \label{PrikryClosed}
$\ordered{\PE, \leq^*}$ is $\gk$-closed.
\end{proposition}
\begin{proof}
Let $\gl < \gk$ and $\ordof {p^\gx} {\gx < \gl} \subseteq \PE$ such that
$\gx_1 > \gx_2 \implies$ $p^{\gx_1} \leq^* p^{\gx_2}$.

By the definition of $\leq^*$ we have
\begin{align*}
\forall \gx_1,\gx_2 < \gl\ 
        \forall \gb \in \supp p^{\gx_1} \intersect \supp p^{\gx_2}\ 
		f^{p^{\gx_1}}(\gb) = f^{p^{\gx_2}}(\gb)
\end{align*}
hence we can set $f = \union \setof{f^{p^{\gx}}} {\gx < \gl}$.

Choose $\ga \gtE \mc(p^\gx)$ for all $\gx < \gl$. Then we set
\begin{align*}
& \forall \gx < \gl\ F^{\prime \gx} = F^{p^\gx} \circ \gp_{\ga, mc(p^\gx)},
\\
& \dom F = \bigintersect_{\gx < \gl} \dom F^{\prime \gx},
\\
& \forall \ordered{\gn_0, \dotsc, \gn_k, \gn} \in \dom F\ 
         F_{\ordered{\gn_0 \dotsc, \gn_k}} (\gn)= 
              \bigunion_{\gx < \gl} 
                        F_{\ordered{\gn_0 \dotsc, \gn_k}}^{\prime \gx} (\gn).
\end{align*}
We shrink $\dom F$, using $\ref{CommutShrink}$, $\gl$ times to get
commutativity with respect to each $p^\gx$. Now we set
$p = \ordered{f, \ga, F}$ and we have
        $\forall \gx < \gl\ p \leq^* p^\gx$.
\end{proof}
The above proposition might give the impression that forcing with 
$\ordered{\PE, \leq^*}$ collapses $\gk^+$. However, this does not happen.
In fact we show that forcing with $\ordered{\PE, \leq^*}$ is the same as forcing
with the Cohen forcing for adding $\power{j(\gk)}$ subsets to $\gk$.
\begin{lemma}
Let $G^*$ be $\ordered{\PE, \leq^*}$-generic. Assume $p \in G^*$.
If $q \in \PE$ is such that $f^q = f^p$ then $q \in G^*$.
\end{lemma}
\begin{proof}
Let $G^*$ be $\ordered{\PE, \leq^*}$-generic,  $p \in G^*$,
$q \in \PE$ is such that $f^q = f^p$. We show that
	$D^* = \setof {r \in \PE} {r \leq^* q}$
is $\leq^*$-dense below p.

Let $s\leq^* p$. Let $\ga \gtE \mc(s), \mc(q)$. By \ref{AnyMC}
there is $r' \leq^* s$ such that $\mc(r') = \ga$.
Since $f^p = f^q$ there is an $E(\ga)$-tree, T, such that
$F^p \circ \gp_{\ga,\mc(p)} \restricted T=
	F^q \circ \gp_{\ga,\mc(q)} \restricted T$.
We set $r$ to be $r'$ with $\dom r = \dom r' \restricted T$.
This insures us that $r \leq^* q$. Hence $D^*$ is dense open below p.

So, there is $r \in G^* \intersect D^*$. Since $r \leq^* q$ we get
$q \in G^*$.
\end{proof}
The above lemma means that the order $\leq^*$ does not separate
between conditions $p, q$ if $f^p=f^q$. Hence we define
\begin{definition}
Let $\PES = \setof {f} {\exists \ga,F\ \ordered{f,\ga,F} \in \PE}$.
We supply $\PES$ with the partial order: $f \leq^* g \iff f \supseteq g$.
\end{definition}
Thus we have the forcing notions $\ordered{\PE, \le}$ (when writing
$\PE$ we mean this forcing), $\ordered{\PE, \leq^*}$, and
$\ordered{\PES, \leq^*}$ (when writing $\PES$ we mean this forcing).
Note that $\PES$ is the Cohen forcing for adding $\power{j(\gk)}$ subsets
to $\gk^+$.
The previous lemma means
\begin{corollary}
Forcing with $\ordered{\PE, \leq^*}$ is the same as forcing with
$\PES$.
\end{corollary}
Note that in \cite{PrikryExtender} the direct extension was not
the Cohen forcing, and was not $\gk^+$-closed. This higher closedness
is used by us in the proofs of \ref{Dichotomy-n}, \ref{NPgeneric}.
\begin{lemma} \label{InducedDensity}
Assume $p \in \PE$, $D$ is $\leq^*$-dense open below 
	$p_{\ordered{\gn_0, \dotsc, \gn_{n-1}}}$. Then
$\setof {f^p \union (f^q \restricted (\dom f^q \setminus \dom f^p)} {q \in D}$
is $\leq^*$-dense open below $f^p$.
\end{lemma}
\begin{proof}
Assume $p \in \PE$, $D$ is $\leq^*$-dense open below 
	$p_{\ordered{\gn_0, \dotsc, \gn_{n-1}}}$. Set
$D^* = \setof {f^p \union (f^q \restricted (\dom f^q \setminus \dom f^p)} {q \in D}$.
We show that $D^*$ is dense open below $f^p$.
Let $f \leq^* f^p$. 

Choose $\ga \gtE \gb$ for all $\gb \in \dom f$.
Set $p^1 = p_{\ordered{\gn_0, \dotsc, \gn_{n-1}}}$.
We enlarge
$F^{p^1}$ so it will be legal to use it in a condition with support
$\supp f$ (as usual this might mean shrinkage of $\dom F^{p^1}$) by choosing
a function $h$ satisfying $j(h)(\ga) = j''(\supp f \setminus \supp f^p)$
and setting
\begin{multline*}
\forall \ordered{\gm_0, \dotsc, \gm_{k-1}} \in \dom F\ 
	F(\gm_0, \dotsc, \gm_{k-1}) = 
\\
	F^{p^1} \circ \gp_{\ga, \mc(p)} 
		(\gm_0, \dotsc, \gm_{k-1})\union h(\gm_{k-1}),
\end{multline*}
\begin{align*}
q = \ordered{f^{p^1} \union f \restricted (\dom f \setminus \dom f^{p^1}), 
	\ga, 
	F
	}.
\end{align*}
Then $q \leq^* p_{\ordered{\gn_0, \dotsc, \gn_{n-1}}}$.
Hence there is $r \leq^* q$, $r \in D$.
This yields $f^p \union (f^r \restricted (\dom f^r \setminus \dom f^p)) \in
	D^*$.
Observing that $f^r \leq^* f^q$ implies 
$f^r \restricted (\dom f^r \setminus \dom f^p) \leq^* f^q \restricted
			(\dom f^q \setminus \dom f^p)$, 
(and $f^q \restricted \dom f \setminus \dom f^p= 
	f \restricted \dom f \setminus \dom f^p$!)
we get 
$f^p \union (f^r \restricted (\dom f^r \setminus \dom f^p)) \leq^* 
	f^p \union (f \restricted (\dom f \setminus \dom f^p)) = f$,
thus proving density.
\end{proof}
\begin{lemma} \label{Dichotomy-n}
Let $D \subseteq \PE$ be dense open, $p \in \PE$, and $n < \gw$. 
Then there is $p^* \leq^* p$, such that
either
\begin{align*}
& \forall \ordered{\gn_0, \dotsc,\gn_{n-1}} \in \dom p^*\ 
         p^*_{\ordered{\gn_0,\dotsc, \gn_{n-1}}}
                         \in D
\intertext{or}
& \forall \ordered{\gn_0, \dotsc,\gn_{n-1}} \in \dom p^*\ 
	\forall q \leq^* p^*_{\ordered{\gn_0, \dotsc, \gn_{n-1}}}\ 
                         q \notin D
\end{align*}
\end{lemma}
\begin{proof}
Assume $D$ is a dense open subset of $\PE$, $p \in \PE$, $n < \gw$.

For each $\ordered{\gn_0, \dotsc, \gn_{n-1}} \in \dom p$ set
\begin{align*}
& D^\in (\gn_0, \dotsc, \gn_{n-1})=\setof 
	{q \leq^* p_{\ordered{\gn_0, \dotsc, \gn_{n-1}}}} 
	{q \in D},
\\
& D^\incompatible(\gn_0, \dotsc, \gn_{n-1}) =\setof
	{r \leq^* p_{\ordered{\gn_0, \dotsc, \gn_{n-1}}}} 
	{\forall q \in D^\in(\gn_0, \dotsc, \gn_{n-1})\  r \incompatible^* q},
\\
& D(\gn_0, \dotsc, \gn_{n-1}) =	D^\in(\gn_0, \dotsc, \gn_{n-1}) \union 
	D^\incompatible(\gn_0, \dotsc, \gn_{n-1}).
\end{align*}
Taking the ultrapower we have
\begin{align*}
&
\begin{aligned}
 D_n^\in= 
	 j_n(D^\in) (\mc(p), & \dotsc, j_{n-1}(\mc(p))) 
	\\
	 & (= \setof 
	{q \leq^* j_n(p)_{\ordered{\mc(p), \dotsc, j_{n-1}(\mc(p))}}} 
	{q \in j_n(D)}),
\end{aligned}
\\
&
\begin{aligned}
 D_n^\incompatible =
  	j_n(D^\incompatible)(\mc(p), & \dotsc,  j_{n-1}(\mc(p))) 
	\\
	& (=\setof
	{r \leq^* j_n(p)_{\ordered{\mc(p), \dotsc, j_{n-1}(\mc(p))}}} 
	{\forall q \in D_n^\in\  r \incompatible^* q}),
\end{aligned}
\\
& D_n =	D_n^\in \union
	D_n^\incompatible.
\end{align*}
Note that the openness of $D$ guarantees us the $\leq^*$-openness of 
$D^\in(\gn_0, \dotsc, \gn_{n-1})$
and that $D^\incompatible(\gn_0, \dotsc, \gn_{n-1})$ is $\leq^*$-open by its definition.
$D(\gn_0, \dotsc, \gn_{n-1})$, as a union of two open sets, is open, and in fact
it is also $\leq^*$-dense
below 
	$p_{\ordered{\gn_0, \dotsc, \gn_{n-1}}}$.
Set
\begin{align*}
& D^*(\gn_0, \dotsc, \gn_{n-1}) =	
\setof {f^p \union (f^q \restricted (\supp q \setminus \supp p))} {q \in D(\gn_0, \dotsc, \gn_{n-1})}.
\\
& 
\begin{aligned}
 D_n^* = j_n(D^*)(\mc(p), & \dotsc, j_{n-1}(\mc(p))) =	
	\\
	& (=\setof {f^{j_n(p)} \union (f^q \restricted 
	(\supp q \setminus \supp j_n(p)))} 
		{q \in D_n)}).
\end{aligned}
\end{align*}
By \ref{InducedDensity}, $D^*(\gn_0, \dotsc, \gn_{n-1})$ is $\leq^*$-dense
open below $f^p$.
Since $\PES$ is $\gk^+$-closed the set 
$D^* = \bigintersect_{\ordered{\gn_0, \dotsc, \gn_{n-1}} \in \dom p}
		D^*(\gn_0, \dotsc, \gn_{n-1})$ is dense open below
$f^p$. Pick $f \in D^*$. We will construct a direct extension of $p$
with $f$ its Cohen part.

Since $f \in D^*$, we have that
	$j_n(f) \in D^*_n$. That is 
there is $q \in D_n$ such that
\begin{align*}	
 &\supp q = \dom(j_n(f)),
 \\
 & f^{q} \restricted \supp q \setminus \supp j_n(p) = j_n(f)
		\restricted \supp q \setminus \supp j_n(p),
 \\
 & f^{q}  \restricted \supp j_n(p) = 
	f^{j_n(p)_{\ordered{\mc(p), \dotsc, 
		j_{n-1}(\mc(p))}}}.
\end{align*}	
We take $\ga \in \dom E$ such that $j_n(\ga) \gtE \mc(q)$, and
there is a function $F'$ such that $j_n(F')(\ga, \dotsc, j_{n-1}(\ga)) = F^q$.
Assume $\mc(q) = j_n(\gb)(\ga, \dotsc, j_{n-1}(\ga))$,
we define a function $F$ with domain an $E(\ga)$-tree so as to have
\begin{align*}
& F = F^p \circ \gp_{\ga, \mc(p)} \restricted \dom F,
\\
& F_{\ordered{\gn_0, \dotsc, \gn_{n-1}}} =
	F'(\gn_0, \dotsc, \gn_{n-1})  \circ
		\gp_{\ga, \gb(\gn_0, \dotsc, \gn_{n-1})}
		\restricted 
 		\dom F_{\ordered{\gn_0, \dotsc, \gn_{n-1}}}.
\end{align*}
We set $p^* = \ordered{f, \ga, F}$. By the construction we have
%\begin{align*}
$j_n(p^*)_{\ordered{\mc(p^*), \dotsc, j_{n-1}(\mc(p^*))}} \leq^* q$,
%\end{align*}
hence
%\begin{align*}
$j_n(p^*)_{\ordered{\mc(p^*), \dotsc, j_{n-1}(\mc(p^*))}} \in D_n$.
%\end{align*}
Since $D_n$ is the union of the two disjoint sets $D^\in_n$, $D^\incompatible_n$,
we can shrink $\dom p^*$ so as to get, after reflecting to $V$, either
\begin{align*}
\forall \ordered{\gn_0, \dotsc,\gn_{n-1}} \in \dom p^*\ 
         p^*_{\ordered{\gn_0,\dotsc, \gn_{n-1}}}
                         \in D^\in(\gn_0, \dotsc, \gn_{n-1})
\end{align*}
or
\begin{align*}
\forall \ordered{\gn_0, \dotsc,\gn_{n-1}} \in \dom p^*\ 
	p^*_{\ordered{\gn_0, \dotsc, \gn_{n-1}}}\ 
                         \in D^\incompatible(\gn_0, \dotsc, \gn_{n-1}).
\end{align*}
Looking at the definition of $D^\incompatible$, $D^\in$ we see that we have
either
\begin{align*}
& \forall \ordered{\gn_0, \dotsc,\gn_{n-1}} \in \dom p^*\ 
         p^*_{\ordered{\gn_0,\dotsc, \gn_{n-1}}}
                         \in D
\intertext{or}
& \forall \ordered{\gn_0, \dotsc,\gn_{n-1}} \in \dom p^*\ 
	\forall q \leq^*p^*_{\ordered{\gn_0, \dotsc, \gn_{n-1}}}\ 
                         q \notin D.
\end{align*}
\end{proof}
\begin{theorem} \label{DenseHomogen}
Let $D \subseteq \PE$ be dense open and $p \in \PE$. 
Then there are $p^* \leq^* p$, $n < \gw$, such that
$\forall \ordered{\gn_0, \dotsc,\gn_{n-1}} \in \dom p^*$
         $p^*_{\ordered{\gn_0,\dotsc, \gn_{n-1}}}
                         \in D$.
\end{theorem}
\begin{proof}
Assume $D$ is a dense open subset of $\PE$, and $p \in \PE$.
Set $p^0 = p$. By induction we construct $p^{n+1}$ by invoking
\ref{Dichotomy-n} for $p^n$ and $n$.
When the induction terminates we have $\ordof {p^n} {n <\gw}$,
a $\leq^*$-decreasing sequence so that
either
\begin{align*}
& \forall \ordered{\gn_0, \dotsc,\gn_{n-1}} \in \dom p^{n+1}\ 
         p^{n+1}_{\ordered{\gn_0,\dotsc, \gn_{n-1}}}
                         \in D
\intertext{or}
& \forall \ordered{\gn_0, \dotsc,\gn_{n-1}} \in \dom p^{n+1}\ 
	\forall q \leq^* p^{n+1}_{\ordered{\gn_0, \dotsc, \gn_{n-1}}}\ 
                         q \notin D.
\end{align*}
Let us pick $p^* \in \PE$ so that $\forall n < \gw\ p \leq^* p^n$.
Noting the openness of $D$ together with
\begin{align*}
\forall \ordered{\gn_0, \dotsc,\gn_{n-1}} \in \dom p^{*}\ 
         p^{*}_{\ordered{\gn_0,\dotsc, \gn_{n-1}}} \leq^*
         p^{n+1}_{\ordered{\gp_{\mc(p^*),\mc(p)}(\gn_0),
			\dotsc, 
			\gp_{\mc(p^*),\mc(p)}(\gn_{n-1})}}
\end{align*}
and the fact that
\begin{multline*}
\forall \ordered{\gn_0, \dotsc,\gn_{n-1}} \in \dom p^{*}\ 
	q \leq^* p^{*}_{\ordered{\gn_0, \dotsc, \gn_{n-1}}} \implies
\\
	q \leq^* p^{n+1}_{\ordered{\gp_{\mc(p^*),\mc(p)}(\gn_0),
			\dotsc, 
			\gp_{\mc(p^*),\mc(p)}(\gn_{n-1})}}
\end{multline*}
we get 
either
\begin{align*}
& \forall n <\gw\ 
\forall \ordered{\gn_0, \dotsc,\gn_{n-1}} \in \dom p^{*}\ 
         p^{*}_{\ordered{\gn_0,\dotsc, \gn_{n-1}}}
                         \in D
\intertext{or}
&\forall n <\gw\ 
\forall \ordered{\gn_0, \dotsc,\gn_{n-1}} \in \dom p^{*}\ 
	\forall q \leq^* p^{*}_{\ordered{\gn_0, \dotsc, \gn_{n-1}}}\ 
                         q \notin D.
\end{align*}
Let us assume, by contradiction that for each $n < \gw$
\begin{align*}
\forall \ordered{\gn_0, \dotsc,\gn_{n-1}} \in \dom p^{*}\ 
	\forall q \leq^* p^{*}_{\ordered{\gn_0, \dotsc, \gn_{n-1}}}\ 
                         q \notin D.
\end{align*}
This is just a cumbersome way to write $\forall q \leq p^*\ q \notin D$,
contradicting the density of $D$.
Hence there is $n < \gw$ such that
\begin{align*}
\forall \ordered{\gn_0, \dotsc,\gn_{n-1}} \in \dom p^{*}\ 
         p^{*}_{\ordered{\gn_0,\dotsc, \gn_{n-1}}}
                         \in D.
\end{align*}
\end{proof}
\begin{claim}
Let $\gs$ be a statement in the $\PE$-forcing language, $p \in \PE$.
Then there is $p^* \leq^* p$ such that $p^* \decides \gs$.
\end{claim}
\begin{proof}
Let $D = \setof {q \in \PE} {q \decides \gs}$. $D$ is a dense open subset of
$\PE$. By \ref{DenseHomogen} there are $p^{*\prime} \leq p$ and $k < \gw$
such that
\begin{align*}
\forall \ordered{\gn_0, \dotsc,\gn_k} \in \dom p^{*\prime}\ 
         p^{*\prime}_{\ordered{\gn_0,\dotsc, \gn_k}}
                         \in D.
\end{align*}
That is
\begin{align*}
\forall \ordered{\gn_0, \dotsc,\gn_k} \in \dom p^{*\prime}\ 
         p^{*\prime}_{\ordered{\gn_0,\dotsc, \gn_k}}
                         \decides \gs.
\end{align*}
Let
\begin{align*}
& A_1 = \setof {\ordered{\gn_0, \dotsc, \gn_k} \in \dom p^{*\prime}}
               {  p^{*\prime}_{\ordered{\gn_0,\dotsc, \gn_k}}
                         \forces \gs},
\\
& A_2 = \setof {\ordered{\gn_0, \dotsc, \gn_k} \in \dom p^{*\prime}}
               {  p^{*\prime}_{\ordered{\gn_0,\dotsc, \gn_k}}
                         \forces \lnot\gs}.
\end{align*}
Obviously, $A_1 \intersect A_2 = \emptyset$. Let $i \in \set{1,2}$
so that $A_i \in E^k(\mc(p^{*\prime}))$.

Let $p^*$ be $p^{*\prime}$ with $\dom p^*$ shrunken to be 
$\dom p^{*\prime} \restricted A_i$.
Then $p^* \decides \gs$.
\end{proof}
So now we know also that in a $\PE$-generic extension
\begin{enumerate}
\item
        There are no new bounded subsets of $\gk$.
\item
        (Hence) No cardinal below $\gk$ is collapsed.
\item
        (Hence) $\gk$ is not collapsed.
\end{enumerate}
We give a direct proof of $\gk^+$ being preserved. Later on
we prove a form of properness of $\PE$ which implies $\gk^+$ is preserved.
\begin{proposition}
$\gk^+$ is preserved in a $\PE$-generic extension.
\end{proposition}
\begin{proof}
Let $\gl < \gk$ and $p \forces \formula {\GN{f} \func \gl \to (\gk^+)^V}$.
For each $\gx < \gl$ let 
        $D_\gx = \setof {q \leq p} {\exists \gz\ q \forces \formula {
                \GN{f}(\VN{\gx}) = \VN{\gz}}}$.
Each $D_\gx$ is a dense open subset of $\PE$ below $p$.
We construct by induction a $\leq^*$-decreasing sequence
        $\ordof {p^\gx} {\gx \leq \gl}$ as follows:
\begin{itemize}
\item
        $p^0 = p$.
\item
        Construct $p^{\gx+1}$, $k^\gx$ from $p^\gx$, $D_\gx$ using 
        \ref{DenseHomogen}.
        Let $g^{\gx} \func [\dom p^{\gx+1}] \restricted [\gk]^{k^\gx} \to 
                        \gk^+$ be defined so that
        $\forall \ordered{\gn_0, \dotsc, \gn_k} \in \dom p^{\gx+1}$
        $p^{\gx+1}_{\ordered{\gn_0,\dotsc, \gn_k{^{\gx}}}} \forces
                \formula{\GN{f}(\gx) = g^\gx(\gn_0, \dotsc, \gn_{k^\gx})
                }$.
\item
        For $\gx$ limit we choose $p^\gx$ such that $p^\gx \leq^* p^{\gx'}$
        for all $\gx' < \gx$.
\end{itemize}
Let $g(\gx, \gn_0,\dotsc, \gn_{k^\gx})=
        g^\gx(\gn_0, \dotsc, \gn_{k^\gx})$. Then
\begin{align*}
\forall {\ordered{\gn_0, \dotsc, \gn_{k^\gx}}} \in \dom p^\gl\ 
        p^\gl_{\ordered{\gn_0, \dotsc, \gn_{k^\gx}}} \forces \formula {
        \GN{f}(\gx) = g(\gx, \gn_0, \dotsc, \gn_{k^\gx})}.
\end{align*}
We set a bound 
\begin{align*}
\gm = \sup \setof {g(\gx, \gn_0, \dotsc, \gn_{k^\gx})} 
                {\gx <\gl,\ \ordered{\gn_0, \dotsc, \gn_{k^\gx}}
                        \in \dom p^\gl}.
\end{align*}
It is clear that $\gm < \gk^+$,
and $p^\gl \forces \formula {\sup \setof {\GN{f}(\gx)} {\gx < \gl} < \VN{\gm}}$.
Hence $\forces_{\PE} \formula{\cf (\gk^+)^V \geq \gk}$.

Since $\forces_{\PE} \formula{ \gk \text{ is singular}}$, we get
$\forces_{\PE} \formula{\cf (\gk^+)^V = (\gk^+)^V}$.
That is $\gk^+$ remains a cardinal in a $\PE$-generic extension.
\end{proof}
All in all we get Gitik-Magidor theorem directly from arbitrary extender:
\begin{theorem} \label{PrikryExtenderMain}
Assume $j \func V \to M \supset M^\gk$, $\crit(j) = \gk$,
$\gen{j} \subset j(\gk)$. Let $E$ be the extender derived from $j$,
and let $G$ be $\PE$-generic. Then in
$V[G]$:
\begin{enumerate}
\item
        $\cf \gk = \gw$.
\item
        $2^\gk = \power{j(\gk)}$.
\item
        No new bounded subsets are added to $\gk$.
\item
        All   cardinals are preserved.
\end{enumerate}
\end{theorem}
For the sake of completeness we show that $\PE$ satisfies a form of
properness. Originally we used this property when we had worked
on the `Radin on Extenders' forcing. In Radin case, $\gk$ might
remain regular and the proof above, of $\gk^+$ being preserved,
fails, so we had to use some other way, thus we stepped
on properness.
Woodin initiated the
use of proper forcing arguments to show cardinal preservation in Radin-style
forcing.

From this point on there are occasional  mentions of $\gc$ which is
`large enough'. As usual this means that whatever we are interested with
appear in $H_\gc$. We also take $N \subelem H_\gc$, $\power{N} = \gk$,
$N \supseteq N^{\upto \gk}$, $\PE \in N$. Note this implies $N \supset \gk+1$.

The notions $\ordered{N,P}$-generic and properness, as defined by
Shelah \cite{ProperForcing}, are as follows:
\begin{definition}
Let $N \subelem H_\gc$ such that
$\power{N}=\gw$,
$P \in N$.
Then $p \in P$ is called $\ordered{N,P}$-generic if
   $p \forces \formula{\forall D\in \VN{N}\ D \text{ is dense open in }\VN{P}\ 
        \implies D \intersect \CN{G} \intersect 
                \VN{N} \not= \emptyset}$.
\end{definition}
\begin{definition}
A forcing notion $P$ is called proper if for all
$N \subelem H_\gc$, $q \in P \intersect N$ such that
        $\power{N}=\gw$,
        $P \in N$,
there is $p \leq q$ which
is $\ordered{N, P}$-generic.
\end{definition}
We adapt these definitions for our needs (namely, larger submodels),
keeping the original names:
\begin{definition}
\label{NewStuff}
Let $N \subelem H_\gc$ such that
        $\power{N}=\gk$,
        $N \supseteq N^{\upto\gk}$,
        $P \in N$.
Then $p \in P$ is called $\ordered{N,P}$-generic if
   $p \forces \formula{\forall D\in \VN{N}\ D \text{ is dense open in }\VN{P}\ 
        \implies D \intersect \CN{G} \intersect 
                \VN{N} \not= \emptyset}$.
\end{definition}
\begin{definition}
A forcing notion $P$ is called proper if for each
$N \subelem H_\gc$, $q \in P \intersect N$ such that
        $\power{N}=\gk$,
        $N \supseteq N^{\upto\gk}$,
        $P \in N$,
there is $p \leq q$ which
is $\ordered{N, P}$-generic.
\end{definition}
\begin{claim} \label{NPgeneric}
Let $p \in \PE$, $N \subelem H_\gc$ such that
$\power{N}=\gk$,
        $N \supseteq N^{\upto\gk}$,
        $p, \PE \in N$.
Then there is $p^* \leq^* p$ such that $p^*$ is $\ordered{N, \PE}$-generic.
\end{claim}
\begin{proof}
Assume $p \in \PE$, $N \subelem H_\gc$ such that
$\power{N}=\gk$,
        $N \supseteq N^{\upto\gk}$,
        $p, \PE \in N$.

Since $\PES$ is $\gk^+$-closed, there is $f \in \PES$ such that for each
$D^* \in N$, a $\leq^*$-dense open subset of $\PES$ below $f^p$, there is 
$g \geq^* f$ such that $g \in D^* \intersect N$. 
We pick
$p^* \leq^* p$ such that $f^{p^*} = f$. We show that $p^*$ is
$\ordered{N, \PE}$-generic.
So let $D \in N$ be dense open subset of $\PE$
and  $q \leq p^*$. 

Then there is
$\ordered{\gn_0, \dotsc, \gn_{n-1}} \in \dom p$ such that 
	$q \leq^*   p_{\ordered{\gn_0, \dotsc, \gn_{n-1}}}$.
We set $D^* = \setof {f^p \union f^r \restricted \supp r \setminus \supp p} {
		r \leq^* p_{\ordered{\gn_0, \dotsc, \gn_{n-1}}},\ 
		r \in D}$.
We note that $D^*$ is $\leq^*$-dense open below $f^p$.
Since $D \in N$, we have $D^* \in N$. By the way we chose $f$ we see that
there is $g \geq^* f$ such that $g \in D^* \intersect N$. Hence
there is $r \in D \intersect N$ such that 
$r \leq^* p_{\ordered{\gn_0, \dotsc, \gn_{n-1}}}$, 
	$f^p \union f^r\restricted
	\supp r \setminus \supp p = g$. In fact $q \compatible^* r$. 
That is there is $q^* \leq^* q$ (a shrinkage of $\dom q$ is enough,
actually) such that $q^* \leq^* r$. Since $r \in D \intersect N$
we get $q^* \forces_{\PE} \formula{\VN{D} \intersect \VN{N}
	\intersect \CN{G} \neq \emptyset}$.
\end{proof}
\begin{corollary} \label{properness}
$\PE$ is proper.
\end{corollary}
\section{Application to $\pcf$ theory.} \label{pcfApplication}
Throughout this section $D$ will be the cofinite filter over $\gw$.
\begin{lemma} \label{GnameToFunc}
If $\gt \in \dom E$, $p \in \PE$, then
there are $p^* \leq^* p$, $n_0 < \gw$, such that
\begin{align*}
\forall n < \gw\ 
j_{\gw}(p^*)_{\ordered{\mc(p^*), \dotsc, j_{n}(\mc(p^*))}} 
                \forces_{j_\gw(\PE)}
         \formula
        {
                j_{\gw}(\CN{G}^\gt)(n_0  + n) = j_{n}(\gt)
        }.
\end{align*}
\end{lemma}
\begin{proof}
Let $\gt \in \dom E$, $p \in \PE$.
By \ref{ArbitraryInSupport} there is $p^* \leq^* p$ 
such that $\gt \in \supp p^*$.
Set $n_0 = \power{f^{p^*}(\gt)}$.
We shrink $\dom p^*$ so that $\forall \ordered{\gn} \in \dom p^*$
$\gt \in F^{p^*}(\gn)$.
Hence, from the definition of $\PE$, we get
$\forall \ordered{\gn_0, \dotsc, \gn_n} \in \dom p^*$
\begin{align*}
p^*_{\ordered{\gn_0, \dotsc, \gn_{n}}} \forces_{\PE} \formula {
        \CN{G}^\gt(n_0 + n) = \gp_{\mc(p^*), \gt} (\gn_{n})
}.
\end{align*}
Translating to ultrapower we get
\begin{align*}
\forall n < \gw\ 
j_{\gw}(p^*)_{\ordered{\mc(p^*), \dotsc, j_{n}(\mc(p^*))}}
        \forces_{j_\gw(\PE)} 
        \formula {j_{\gw}(\CN{G}^\gt)(n_0 + n) = j_{n}(\gt)}.
\end{align*}
\end{proof}
We would have liked to have $\gr < \gt \implies G^\gr/D < G^\gt/D$.
However, the Cohen start-segments of $G^\gr$, $G^\gt$ ruin this.
We can get a good approximation to this monotonicity using shifts
of $G^\gr$, hence the following definition.
By $\sZ$ we mean the set of integers
        $\set{0, 1, -1, 2, -2, \dotsc}$.
\begin{definition}
Assume $\gt \in \dom E$, $k \in \sZ$. Then, in $V[G]$,
 $G^{\gt,k}\func \gw \to \gk$ is 
\begin{align*}
G^{\gt,k}(n) =
         \begin{cases}
        G^\gt(n - k)    & k \leq n < \gw,
        \\
        0               & 0 \leq n < k.
        \end{cases}
\end{align*}
As usual, $\CN{G}^{\gt,k}$ will be the $\PE$-name of this function.
\end{definition}
\begin{lemma}
$\forall \gt \in \dom E$ $\forall k_1, k_2 \in \sZ$ 
        $\cf \prod G^{\gt,k_1}/D = \cf \prod G^{\gt,k_2}/D$.
\end{lemma}
\begin{proof}
\end{proof}
\begin{lemma}
If $\gt \in \dom E$, $k_1, k_2 \in \sZ$, $k_1 < k_2$, 
        then $G^{\gt,k_1}/D > G^{\gt, k_2}/D$.
\end{lemma}
\begin{proof}
This is immediate since $G^\gt$ is a strictly increasing sequence.
\end{proof}
\begin{lemma} \label{CannonicalK<}
If $p\in \PE$, $\gr,\gt \in \dom E$, $\gr<\gt$, then
there are $p^* \leq^* p$, $k \in \sZ$, such that
$p^* \forces_{\PE} \formula{\CN{G}^{\gr,k}/D < \CN{G}^\gt/D <
        \CN{G}^{\gr,k-1}/D}$.
\end{lemma}
\begin{proof}
Let $\gr, \gt \in \dom E$, $\gr < \gt$, $p \in \PE$.
By \ref{ArbitraryInSupport}, \ref{GnameToFunc}, there are 
$p^* \leq^* p$, $n_0, n_1 < \gw$, such that
$\gr, \gt \in \supp p^*$, and
\begin{align*}
& \forall n < \gw\ j_\gw(p^*)_{\ordered{\mc(p^*), \dotsc, j_n(p^*)}}
                \forces_{j_\gw(\PE)}
        \formula{j_\gw(G^\gr)(n_0 + n) = j_n(\gr)},
\\
& \forall n < \gw\ j_\gw(p^*)_{\ordered{\mc(p^*), \dotsc, j_n(p^*)}}
                \forces_{j_\gw(\PE)}
        \formula{j_\gw(G^\gt)(n_1 + n) = j_n(\gt)}.
\end{align*}
We shrink $\dom p^*$ in order to have
$\forall \ordered{\gn} \in \dom p^*$
         $\gr \in F^{p^*}(\gn) \iff \gt \in F^{p^*}(\gn)$.
We set $k=n_1 - n_0$.
Hence
\begin{multline*}                                                
\forall n < \gw\ 
j_\gw(p^*)_{\ordered{\mc(p^*), \dotsc, j_{n+1}(p^*)}}
        \forces_{j_\gw(\PE)}
        \formula{j_\gw(\CN{G}^{\gr,k})(n_1 + n) = 
\\
                j_\gw(\CN{G}^{\gr})(n_1 + n - (n_1 - n_0))=
                j_\gw(\CN{G}^{\gr})(n_0 + n)=
\\
                j_{n}(\gr) < 
                j_n(\gt) = j_\gw(\CN{G}^\gt)(n_1 + n)
},
\end{multline*}
and 
\begin{multline*}                                                
\forall n < \gw\ 
j_\gw(p^*)_{\ordered{\mc(p^*), \dotsc, j_{n+1}(p^*)}}
        \forces_{j_\gw(\PE)}
        \formula{
                j_\gw(\CN{G}^\gt)(n_1 + n) = 
\\
                j_{n}(\gt) <
                j_{n+1}(\gr) = 
                j_\gw(\CN{G}^\gr)(n_0 + n + 1) =
\\
                j_\gw(\CN{G}^{\gr})(n_1 + n - (n_1 - n_0 - 1)=
                j_\gw(\CN{G}^{\gr,k-1})(n_1 + n)
}.
\end{multline*}
Reflection to $V$ and shrinking $\dom p^*$ a bit yield
        $\forall \ordered{\gn_0, \dotsc, \gn_{n+1}} \in \dom p^*$
\begin{align*}
p^*_{\ordered{\gn_0, \dotsc, \gn_{n+1}}}
        \forces_{\PE}
        \formula{\CN{G}^{\gr,k}(n_1 + n)   = 
                & \gp_{\mc(p^*), \gr}(\gn_{n}) <
                \\
                & \gp_{\mc(p^*), \gt}(\gn_n) = \CN{G}^\gt(n_1 + n)
        },
\end{align*}
and
\begin{align*}
p^*_{\ordered{\gn_0, \dotsc, \gn_{n+1}}}
        \forces_{\PE}
        \formula{\CN{G}^{\gt}(n_1 + n)   = 
                & \gp_{\mc(p^*), \gt}(\gn_{n}) <
                \\
                & \gp_{\mc(p^*), \gr}(\gn_{n+1}) = \CN{G}^{\gr,k-1}(n_1 + n)
        }.
\end{align*}
Which means
        $p^* \forces_{\PE} \formula{
                \CN{G}^{\gr, k}/D < \CN{G}^\gt/D <
                        \CN{G}^{\gr, k-1}/D
                }$.
\end{proof}
The following lemma is the same as the previous one with $\gr > \gt$
substituted for $\gr < \gt$.
\begin{lemma} \label{CannonicalK>}
If $p\in \PE$, $\gr,\gt \in \dom E$, $\gr > \gt$, then
there are $p^* \leq^* p$, $k \in \sZ$, such that
$p^* \forces_{\PE} \formula{\CN{G}^{\gr,k}/D < \CN{G}^\gt/D <
        \CN{G}^{\gr,k-1}/D}$.
\end{lemma}
\begin{proof}
Let $\gr, \gt \in \dom E$, $\gr < \gt$, $p \in \PE$.
By \ref{ArbitraryInSupport}, \ref{GnameToFunc}, there are 
$p^* \leq^* p$, $n_0, n_1 < \gw$, such that
$\gr, \gt \in \supp p^*$, and
\begin{align*}
& \forall n < \gw\ j_\gw(p^*)_{\ordered{\mc(p^*), \dotsc, j_n(p^*)}}
                \forces_{j_\gw(\PE)}
        \formula{j_\gw(G^\gr)(n_0 + n) = j_n(\gr)},
\\
& \forall n < \gw\ j_\gw(p^*)_{\ordered{\mc(p^*), \dotsc, j_n(p^*)}}
                \forces_{j_\gw(\PE)}
        \formula{j_\gw(G^\gt)(n_1 + n) = j_n(\gt)}.
\end{align*}
We shrink $\dom p^*$ in order to have
$\forall \ordered{\gn} \in \dom p^*$
         $\gr \in F^{p^*}(\gn) \iff \gt \in F^{p^*}(\gn)$.
We set $k=n_1 - n_0 + 1$.
Hence
\begin{multline*}                                                
\forall n < \gw\ 
j_\gw(p^*)_{\ordered{\mc(p^*), \dotsc, j_{n+1}(p^*)}}
        \forces_{j_\gw(\PE)}
        \formula{j_\gw(\CN{G}^{\gr,k})(n_1 + n + 1) = 
\\
                j_\gw(\CN{G}^{\gr})(n_1 + n  + 1 - (n_1 - n_0 + 1))=
                j_\gw(\CN{G}^{\gr})(n_0 + n)=
\\
                j_{n}(\gr) < 
                j_{n+1}(\gt) = j_\gw(\CN{G}^\gt)(n_1 + n + 1)
},
\end{multline*}
and 
\begin{multline*}                                                
\forall n < \gw\ 
j_\gw(p^*)_{\ordered{\mc(p^*), \dotsc, j_{n+1}(p^*)}}
        \forces_{j_\gw(\PE)}
        \formula{
                j_\gw(\CN{G}^\gt)(n_1 + n + 1) = 
\\
                j_{n+1}(\gt) <
                j_{n+1}(\gr) = 
                j_\gw(\CN{G}^\gr)(n_0 + n + 1) =
\\
                j_\gw(\CN{G}^{\gr})(n_1 + n + 1 - (n_1 - n_0)=
                j_\gw(\CN{G}^{\gr,k-1})(n_1 + n + 1)
}.
\end{multline*}
Reflection to $V$ and shrinking $\dom p^*$ a bit yield
        $\forall \ordered{\gn_0, \dotsc, \gn_{n+1}} \in \dom p^*$
\begin{align*}
p^*_{\ordered{\gn_0, \dotsc, \gn_{n+1}}}
        \forces_{\PE}
        \formula{\CN{G}^{\gr,k}(n_1 + n + 1)   = 
                & \gp_{\mc(p^*), \gr}(\gn_{n}) <
                \\
                & \gp_{\mc(p^*), \gt}(\gn_{n+1}) = \CN{G}^\gt(n_1 + n + 1)
        },
\end{align*}
and
\begin{align*}
p^*_{\ordered{\gn_0, \dotsc, \gn_{n+1}}}
        \forces_{\PE}
        \formula{\CN{G}^{\gt}(n_1 + n + 1)   = 
                & \gp_{\mc(p^*), \gt}(\gn_{n + 1}) <
                \\
                & \gp_{\mc(p^*), \gr}(\gn_{n+1}) = \CN{G}^{\gr,k-1}(n_1 + n + 1)
        }.
\end{align*}
Which means
        $p^* \forces_{\PE} \formula{
                \CN{G}^{\gr, k}/D < \CN{G}^\gt/D <
                        \CN{G}^{\gr, k-1}/D
                }$.
\end{proof}
\begin{definition}
For each $\gt \in \dom E$ we set $G^{*\gt} = G^{\gt,k}$ where $k \in \sZ$
is chosen so that $G^{\gt,k}/D < G^{\gk, -1}/D < G^{\gt, k-1}$.
Again, $\CN{G}^{*\gt}$ is the $\PE$-name of $G^{*\gt}$.
\end{definition}
Immediate corollary of this definition and \ref{GnameToFunc} is
\begin{corollary} \label{GnameToFuncPar}
Assume $p \in \PE$, $\gr, \gt \in \dom E$ Then there are $p^* \leq^* p$,
$n_0 < \gw$
such that $\gr, \gt \in \supp p^*$ and 
\begin{align*}
& \forall n < \gw \ 
        j_\gw(p^*)_{\ordered{\mc(p), \dotsc, j_n(\mc(p))}} 
                \forces_{j_\gw(\PE)} \formula
        {
                j_\gw(\CN{G}^{*\gr}) = j_n(\gr)
        },
\\
& \forall n < \gw \ 
        j_\gw(p^*)_{\ordered{\mc(p), \dotsc, j_n(\mc(p))}} 
                \forces_{j_\gw(\PE)} \formula
        {
                j_\gw(\CN{G}^{*\gt})= j_n(\gt)
        }.
\end{align*}
\end{corollary}
\begin{corollary} \label{IncreasingSequence}
\ordof {G^{*\gt}/D} {\gt \in \dom E} is a strictly increasing sequence.
\end{corollary}
\begin{proof}
Let $p \in \PE$, $\gr, \gt \in \dom E$, $\gr < \gt$.
By \ref{GnameToFuncPar} 
there are $p^* \leq^* p$, $n_0 < \gw$, such that
$\gr, \gt \in \supp p^*$, and
\begin{align*}
& \forall n < \gw\ j_\gw(p^*)_{\ordered{\mc(p^*), 
                                                        \dotsc, j_n(p^*)}}
                \forces_{j_\gw(\PE)}
        \formula{j_\gw(G^{*\gr})(n_0 + n) = j_{n}(\gr)},
\\
& \forall n < \gw\ j_\gw(p^*)_{\ordered{\mc(p^*), \dotsc, j_n(p^*)}}
                \forces_{j_\gw(\PE)}
        \formula{j_\gw(G^{*\gt})(n_0 + n) = j_{n}(\gt)}.
\end{align*}
Since $\gr < \gt$ we have $\forall n < \gw$ $j_n(\gr) < j_n(\gt)$, hence
\begin{align*}
& \forall n < \gw\ j_\gw(p^*)_{\ordered{\mc(p^*), 
                                                        \dotsc, j_n(p^*)}}
                \forces_{j_\gw(\PE)}
        \formula{j_\gw(G^{*\gr})(n_0 + n) <
                        j_\gw(G^{*\gt})(n_0 + n)}.
\end{align*}
Reflecting to $V$ and shrinking $p^*$ a bit yield
\begin{align*}
& \forall \ordered{\gn_0, \dotsc, \gn_n} \ 
        p^*_{\ordered{\gn_0, \dotsc, \gn_n}}
                \forces_{\PE}
                        \formula{G^{*\gr}(n_0 + n) <
                                G^{*\gt}(n_0 + n)}.
\end{align*}
Hence $p^* \forces_{\PE} \formula{ \CN{G}^{*\gr}/D < \CN{G}^{*\gt}/D}$.
\end{proof}
\begin{lemma} \label{NameToFunc}
If $\gt \in \dom E$,
$p \forces_{\PE} \formula {\GN{f} \in \prod \CN{G}^\gt}$,
then there are $p^* \leq^* p$, $n_0 < \gw$, 
$\ordof{\ga_n} {n < \gw} \in \prod_{n<\gw} j_n(\gt)$,
such that
\begin{align*}
\forall n < \gw\ 
j_{\gw}(p^*)_{\ordered{\mc(p^*), \dotsc, j_{n}(\mc(p^*))}} 
         \forces_{j_\gw(\PE)}
         \formula
        {
                j_{\gw}(\GN{f})(n_0 + n) = \ga_n
        }.
\end{align*}
\end{lemma}
\begin{proof}
Let $\gt \in \dom E$,
$p \forces_{\PE} \formula {\GN{f} \in \prod \CN{G}^\gt}$.
By \ref{GnameToFunc} there are $q \leq^* p$, $n_0 < \gw$, such that
\begin{align} \label{eq:GNameToFunc}
\forall n < \gw\ 
j_{\gw}(q)_{\ordered{\mc(q), \dotsc, j_{n}(\mc(q))}} 
                \forces_{j_\gw(\PE)}
         \formula
        {
                j_{\gw}(\CN{G}^\gt)(n_0  + n) = j_{n}(\gt)
        }.
\end{align}
We construct by induction a $\leq^*$-decreasing sequence 
$\ordof {r^n} {n < \gw}$, and $\ordof {\ga_n} {n < \gw} \in 
                \prod_{n<\gw} j_n(\gt)$ as follows:
\begin{itemize}
\item
        $n = 0$: $r^0 = q$.
\item
        $n= m+1$:
        Let $D_m = \setof {r \in \PE} {\exists \gz<\gk \ r\forces_{\PE} \formula
                                        {\GN{f}(n_0 + m) = \gz}}$.
        By \ref{DenseHomogen} there are $r^{m+1} \leq^* r^m$,
         $k < \gw$, 
        $f\func \dom r^{m+1} \intersect [\gk]^{k} \to \gk$, such that
        $r^{m+1}_{\ordered{\gn_0, \dotsc, \gn_{k}}} \forces_{\PE} \formula {
                \GN{f}(n_0 + m) = f(\gn_0, \dotsc, \gn_{k})
        }$. In ultrapower language we get
        \begin{multline} \label{eq:GeneralRepresen}
        j_\gw(r^{m+1})_{\ordered{\mc(r^{m+1}), \dotsc, j_k(\mc(r^{m+1}))}}
                \forces_{j_\gw(\PE)} 
                \\
                \formula {
                j_{\gw}(\GN{f})(n_0 + m) =
                \\
                j_\gw(f)(\mc(r^{m+1}), \dotsc, j_k(\mc(r^{m+1})))
                }.
        \end{multline}
         From (\ref{eq:GNameToFunc}) we infer
        \begin{multline*}
        j_\gw(r^{m+1})_{\ordered{\mc(r^{m+1}), \dotsc, j_k(\mc(r^{m+1}))}}
                \forces_{j_\gw(\PE)} 
                \\
                \formula {
                j_\gw(f)(\mc(r^{m+1}), \dotsc, j_k(\mc(r^{m+1}))) =
                \\
                j_{\gw}(\GN{f})(n_0 + m) <
                j_{\gw}(\CN{G}^\gt)(n_0 + m) = j_{m}(\gt)
        }.
        \end{multline*}
 This, of course, just means
        \begin{align*}
         j_\gw(f)(\mc(r^{m+1}), \dotsc, j_k(\mc(r^{m+1}))) < 
                j_{m}(\gt).
        \end{align*}
  We set $\ga_m = j_\gw(f)(\mc(r^{m+1}), \dotsc, j_k(\mc(r^{m+1})))$.
        Since $\ga_m < j_{m}(\gt)$,
        there is $g\func \dom r^{m+1} \intersect [\gk]^{m+1} \to \gk$ such that
        \begin{align*}
         \ga_m = j_\gw(g)(\mc(r^{m+1}), \dotsc, j_{m}(\mc(r^{m+1}))).
        \end{align*}
        So in \ref{eq:GeneralRepresen}, we can substitute $g$ for $f$
        and $m$ for $k$ and get
        \begin{align*}
        j_\gw(r^{m+1})_{\ordered{\mc(r^{m+1}), \dotsc, j_m(\mc(r^{m+1}))}}
                \forces_{j_\gw(\PE)} 
                \formula {
                j_{\gw}(\GN{f})(n_0 + m) = \ga_m
                }.
        \end{align*}
\end{itemize}
We get the conclusion by using \ref{PrikryClosed} to find $p^* \in \PE$
such that $\forall n < \gw$ $p^* \leq r^n$. 
\end{proof}
\begin{lemma} \label{OrdinalCf}
Assume  $\gt \in \dom E$, $\cf \gt > \gw$, $\cf \gt \neq \gk$, 
$\ordof {\ga_n} {n < \gw} \in \prod_{n < \gw} j_n(\gt)$.
Then there is $\gr < \gt$ such that
$\ordof {\ga_n} {n < \gw} < \ordof{j_n(\gr)} {n < \gw}$.
\end{lemma}
\begin{proof}
We split the proof according to the relation between $\cf \gt$ and $\gk$:
\begin{itemize}
\item 
        $\cf \gt > \gk$:
        We note that for each $n < \gw$ there is $f_n\func [\gk]^n \to \gt$, 
        $\gb_n \in \dom E$ such that 
        $j^n(f_n)(\gb_n, \dotsc, j_{n-1}(\gb_n)) = \ga_n$.
Since $\cf \gt > \gk$, there is $\gr < \gt$ such that
        $\forall n < \gw\ 
                \forall \ordered{\gn_0, \dotsc, \gn_{n-1}} \in [\gk]^n$
                $\gr > f_n(\gn_0, \dotsc, \gn_{n-1})$.
Hence $\forall n < \gw\ j_n(\gr) > \ga_n$.
\item
        $\cf \gt < \gk$:
        Let $A=\ordof {\gt_\gx} {\gx < \cf \gt}$ be cofinal in $\gt$.
        So for each $n < \gw$ we get $j_n(A)$ is cofinal in $j_n(\gt)$.
         Since $\cf \gt < \gk$ we have that $j_n(A) = j_ni A$.
        This means that for each $n < \gw$ there is $\gx_n < \cf \gt$
        such that $\ga_n < j_n(\gt_{\gx_n}) < j_n(\gt)$.
        Since $\cf \gt > \gw$ there is $\gx < \cf \gt$ such that
        for all $n < \gw$ $\gx > \gx_n$. Let $\gr = \gt_\gx$.
        Then for all $n < \gw$ $\gr > \gt_{\gx_n}$, and
        $j_n(\gr) > j_n(\gt_{\gx_n}) > \ga_n$. Hence
                $\ordof {\ga_n} {n < \gw} < \ordof {j_n(\gr)} {n < \gw}$.
\end{itemize}
\end{proof}
\begin{corollary}
If $\gt \in \dom E$, $\cf \gt > \gw$, $\cf \gt \neq \gk$, then
        $\forces_{\PE} \formula{
                \tcf \prod \CN{G}^{*\gt}/D = \cf \gt
                }$.
\end{corollary}
\begin{proof}
Let $\gt \in \dom E$, $\cf \gt > \gw$, $\cf \gt \neq \gk$.
By \ref{IncreasingSequence}, $\ordof{G^{*\gr}/D} {\gr < \gt}$ is a strictly
increasing sequence below $G^{*\gt}/D$.
We will get the conclusion of the lemma if we prove
\begin{align*}
\forces_{\PE} \formula {
        \GN{f} \in \prod \CN{G}^{*\gt} \implies \exists \gr <  \gt\ 
        \GN{f}/D < \CN{G}^{*\gr} / D
        }.
\end{align*}
So, let $p \forces_{\PE} \formula {\GN{f} \in \prod \CN{G}^{*\gt}}$.

By \ref{NameToFunc}, there are $p^* \leq^* p$, $n_0 < \gw$,
        $\ordof{\ga_n} {n < \gw} \in \prod_{n<\gw} j_n(\gt)$ such that
\begin{align*}
& j_\gw(p^*)_{\ordered{\mc(p^*), \dotsc, j_n(\mc(p^*))}} \forces_{j_\gw(\PE)}
        \formula{
        j_\gw(\GN{f})(n_0 + n) = \ga_n
        },
\\
& j_\gw(p^*)_{\ordered{\mc(p^*), \dotsc, j_n(\mc(p^*))}} \forces_{j_\gw(\PE)}
        \formula{
        j_\gw(\CN{G}^{*\gt})(n_0 + n) = j_n(\gt)
        }.
\end{align*}
By \ref{OrdinalCf}, there is $\gr < \gt$ such that
        $\ordof {\ga_n} {n < \gw} < \ordof {j_n(\gr)} {n < \gw}$.
By \ref{GnameToFuncPar} there is $p^{**} \leq^* p^*$  such that
\begin{align*}
j_\gw(p^{**})_{\ordered{\mc(p^{**}), \dotsc, j_n(\mc(p^{**}))}} 
        \forces_{j_\gw(\PE)}
        \formula{
                j_\gw(\CN{G}^{*\gr})(n_0 + n) = j_n(\gr)
        }.
\end{align*}
Reflecting to $V$ and shrinking $\dom p^{**}$ yield
$p^{**}     \forces_{\PE}
        \formula{
        \GN{f} / D < \CN{G}^{*\gr}/D
        }$.
\end{proof}
The situation shown in the following theorem was suggested to us by M.~Gitik,
it summarizes the facts proved previously for the $E$ is a superstrong
extender case.
\begin{theorem} \label{ScaleTheorem}
Assume GCH,
        $j\func V \to M \supset M^\gk$, $M \supset V_{j(\gk)}$, 
	$\crit{j} = \gk$,
		$\gen{j} \subset j(\gk)$.
Let $E$ by the extender derived from $j$ and
$G$ be $\PE$-generic.
Then $V[G]$ is a cardinal preserving generic extension in which 
$\cf \gk = \gw$,
$\gk^\gw = j(\gk)$, and
$\forall \gl \in [\gk, j(\gk)]_{\text{Reg}}$ there is a function
$G^{\gl}\func \gw \to \gk$ such that
$\ordof {G^{\gl}/D} {\gl \in [\gk,j(\gk))_{\text{Reg}}}$ is a scale in
        $\prod G^{j(\gk)}$,
and $\forall \gl \in [\gk, j(\gk))_{\text{Reg}}$
$\tcf \prod G^{\gl}/D = \gl$.
\end{theorem}
For the sake of completeness we analyze the cofinality of
$\prod G^{*\gt}$ when $\cf \gt = \gw$.
\begin{lemma} \label{OrdinalCfk}
Assume  $\gt \in \dom E$, $\cf \gt = \gk$, 
$\ordof {\ga_n} {n < \gw} \in \prod_{n < \gw} j_n(\gt)$.
Then there is $\gr < j(\gt)$ such that
$\ordof {\ga_n} {n < \gw} /D< 
        \ordered{0} \append \ordof{j_n(\gr)} {0 < n < \gw} /D$.
\end{lemma}
\begin{proof}
We fix $n > 0$. Then there are $f\func [\gk]^n \to \gt$, $\gb \in \dom E$,
such that
\begin{align*}
j_n(f)(\gb, \dotsc, j_{n-1}(\gb)) = \ga_n < j_n(\gt),
\end{align*}
where $\forall \ordered{\gn_0, \dotsc, \gn_{n-2}} \in [\gk]^{n-2}$
        $j(f)(\gn_0, \dotsc, \gn_{n-2}, \gb) < j(\gt)$.
Since $\cf j(\gt) > \gk$, there is $\gr < j(\gt)$ such that
$\forall \ordered{\gn_0, \dotsc, \gn_{n-2}} \in [\gk]^{n-2}$
        $j(f)(\gn_0, \dotsc, \gn_{n-2}, \gb) < \gr$.
Hence
        $\ga_n = j_n(f)(\gb, \dotsc, j_{n-1}(\gb)) < j_{n-1}(\gr) < j_n(\gt)$.

So, for each $n > 0$ there is $\gr_n < j(\gt)$ such that 
$\ga_n < j_{n-1}(\gr_n)$.
Since $\cf j(\gt) > \gk$, there is $\gr < j(\gt)$ such that
        $\forall n > 0$ $\gr > \gr_n$. Hence
        $\forall n > 0$  $\ga_n < j_{n-1}(\gr) < j_n(\gt)$. That is
$\ordof{\ga_n}  {n < \gw} /D < \ordered{0}  \append
                        \ordof {j_{n-1}(\gr)} {0 < n < \gw} /D$.
\end{proof}
\begin{corollary}
If $\gt \in \dom E$, $\cf \gt = \gk$, then
        $\forces_{\PE} \formula{
                \tcf \prod \CN{G}^{*\gt}/D = \cf j(\gk)
                }$.
\end{corollary}
\begin{proof}
Assume $\gt \in \dom E$, $\cf \gt = \gk$.
For each $\gr < j(\gt)$ there are $\gt_\gr \in \dom E$, 
and a function $\Bar{h}_\gr \func \gk \to \gt$
such that $\gr = j(\Bar{h}_\gr)(\gt_\gr)$.
In the generic extension we set 
        $h_\gr = \ordered{0} \union 
                \Bar{h}''_\gr G^{*{\gt_\gr}}$. 
We note that
        $\ordof {h_\gr/D} {\gr < j(\gt)}$ is an increasing
sequence in  $\prod_{n < \gw} j_n(\gt)$. 
We will get the conclusion of the lemma if we prove
$\forces_{\PE} \formula {
        \GN{f} \in \prod \CN{G}^{*\gt} \implies \exists \gr <  j(\gt)\ 
        \GN{f}/D < h_\gr / D
        }$.
So, let $p \forces_{\PE} \formula {\GN{f} \in \prod \CN{G}^{*\gt}}$.

By \ref{NameToFunc}, there are $p^* \leq^* p$, $n_0 < \gw$,
        $\ordof{\ga_n} {n < \gw} \in \prod_{n<\gw} j_n(\gt)$ such that
\begin{align*}
& j_\gw(p^*)_{\ordered{\mc(p^*), \dotsc, j_n(\mc(p^*))}} \forces_{j_\gw(\PE)}
        \formula{
        j_\gw(\GN{f})(n_0 + n) = \ga_n
        },
\\
& j_\gw(p^*)_{\ordered{\mc(p^*), \dotsc, j_n(\mc(p^*))}} \forces_{j_\gw(\PE)}
        \formula{
        j_\gw(\CN{G}^{*\gt})(n_0 + n) = j_n(\gt)
        }.
\end{align*}
By \ref{OrdinalCfk}, there is $\gr < j(\gt)$ such that
        $\ordof {\ga_n} {n < \gw}/D < \ordered{0} \append
                                \ordof{j_{n-1}(\gr)} {n < \gw} /D$.
By \ref{GnameToFuncPar} there is $p^{**} \leq^* p^*$  such that
\begin{align*}
j_\gw(p^{**})_{\ordered{\mc(p^{**}), \dotsc, j_n(\mc(p^{**}))}} 
        \forces_{j_\gw(\PE)}
        \formula{
                j_\gw(\CN{G}^{*\gt_\gr})(n_0 + n) = j_n(\gt_\gr)
        }.
\end{align*}
Hence
\begin{multline*}
j_\gw(p^{**})_{\ordered{\mc(p^{**}), \dotsc, j_n(\mc(p^{**}))}} 
        \forces_{j_\gw(\PE)}
        \formula{
                j_\gw(h_\gr)(n_0 + n) = j_\gw(\Bar{h}_\gr)(j_n(\gt_\gr))=
                \\
                j_{n-1}(\gr)
        }.
\end{multline*}
Reflecting to $V$ and shrinking $\dom p^{**}$ yield
$p^{**}     \forces_{\PE}
        \formula{
        \GN{f} / D < \VN{h}_{\gr}/D
        }$.
\end{proof}
The only thing we are able to say when $\cf \gt =\gw$ is
	$\cf \prod G^{*\gt} \leq 2^{\aleph_0}$.
%
%
% Generic by Iteration
%
\section{Generic by Iteration} \label{GenericByIteration}
What we would have liked to have is
\begin{Expected*}
Assume GCH, $j \func V \to M \supset M^\gk$, $\crit(j) = \gk$,
	$\gen{j} \subset j(\gk)$.
Let $E$ be the extender derived from $j$.
Then there is $G \in V$ which is $j_\gw(\PE)$-generic over
$M_\gw$.
\end{Expected*}
Alas, we were not able to get this theorem. We have two approximations and
one `fake'. Namely, we can either find a generic filter over an 
elementary submodel,
or we can start from a stronger assumption, 
or -- `faking' -- we force over $V$ with Cohen forcing
to get a generic over $M_\gw$.
\begin{theorem}
Assume GCH, $j \func V \to M \supset M^\gk$, $\crit(j) = \gk$,
	$\gen{j} \subset j(\gk)$.
Let $E$ be the extender derived from $j$.
Let $N \in M_\gw$ be such that in $M_\gw$ we have: $N \subelem H^{M_\gw}_\gc$, 
$\power{N} = \gk_\gw$,
$N \supseteq N^{\upto \gk_\gw}$. Then there is $G \in V$ which is $j_\gw(\PE)$-generic
over $N$.
\end{theorem}
\begin{proof}
By \ref{NPgeneric} invoked in $M_\gw$, there is $p \in j_\gw(\PE)$ which forces
genericity over $N$. Let $k < \gw$ be such that $p = j_{k,\gw}(p^k)$,
$N = j_{k,\gw}(N^k)$. Of course for each $k \leq n <\gw$ there are
$p^n$, $N^n$ such that $j(p^n) = p$, $j(N^n)= N$.

We set
\begin{align*}
& G^n = \setof {q \in j_n(\PE)} {q \geq 
                j_{k,n}(p^k)_{\ordered{\mc(p^k), \dotsc, j_{k,n-1}(\mc(p^k))}}},
\\
& G = \bigunion_{k\leq n < \gw} j''_{n,\gw} G^n.
\end{align*}
That $G$ is a filter is seen immediately. We show that it intersects
each dense open subset in $N$. So, let $D \in N$ be a dense open subset
of $j_\gw(\PE)$.

Let $D^n$, $n$ be such that $k \leq n < \gw$, $j_{n,\gw}(D^n) = D$.
Of course, $p^n$ witness genericity over $N^n$.
 That means, there are
$q \geq p^n$, $l < \gw$ such that
$\forall \ordered{\gn_0,\dotsc,\gn_{l-1}} \in \dom q$
                $q_{\ordered{\gn_0, \dotsc, \gn_{l-1}}} \in D^n$.
Hence $j_{n,n+l}(q)_{\ordered{\mc(q), \dotsc, j_{n,n+l-1}(\mc(q))}} \in 
        j_{n,n+l}(D^n)$.
Hence $j_{n,\gw}(q)_{\ordered{\mc(q), \dotsc, j_{n,n+l}(\mc(q))}}
 \in D$.
\end{proof}
\begin{theorem}
Assume GCH,
$i\func V \to N \supset N^\gl$, $\crit(i) = \gk$, 
$\gl = \power{i(\gk)}^{(\gk^+)}$.
Let $F$ be the extender derived from $i$, and
$E = F \restricted i(\gk)$.
Then there is $G \in V$ which is $i_\gw(\PE)$-generic
over $N_\gw$.
\end{theorem}
\begin{proof}
Let $\fA_n = \setof {A \in N_n} {A\  \text{is a maximal anti-chain in}\ i_n(\PE)}$.
For each a maximal anti-chain $A$ we set 
	$D(A) = \setof{a \in A} {p \leq a}$.
The $\gl$-super-compactness means that $i''_{n,n+1} \fA_n \in N_{n+1}$.
Moreover, for $q \in i_{n+1}(\PE)$ we can invoke 
\ref{DenseHomogen} in $N_{n+1}$ for $i_{n}(\gl)$-times to get $p \leq^* q$
such that for each $A \in \fA_n$ there is $l < \gw$ such that
$\forall \ordered{\gn_0, \dotsc, \gn_{l-1}} \in \dom p$
        $p_{\ordered{\gn_0,\dotsc,\gn_{l-1}}} \in i_{n,n+1}(D(A))$.

Using this fact we construct a sequence $\ordof{p^n} {n < \gw}$ such that
$p^0 \in \PE$ is arbitrary, $p^{n+1} \leq^* i_{n,n+1}(p^n)_{\ordered{\mc(p^n)}}$
and for each $A \in \fA_n$ there is $l < \gw$ such that
$\forall \ordered{\gn_0, \dotsc, \gn_{l-1}}\in \dom p^{n+1}$
 $p^{n+1}_{\ordered{\gn_0, \dotsc, \gn_{l-1}}} \in i_{n, n+1}(D(A))$.

Let $G_n = \setof {q \in i_n(\PE)} {q \geq p^n}$,
        $G = \bigunion_{n < \gw} i''_{n,\gw} G_n$. 
We show $G$ is $i_\gw(\PE)$-generic
over $N_\gw$. Let $A \in N_\gw$ be a maximal anti-chain in $i_\gw(\PE)$.

Take $A_n$, $n < \gw$ such that $i_{n,\gw}(A_n) = A$.
Then there is $l < \gw$ such that
$\forall \ordered{\gn_0, \dotsc, \gn_{l-1}}\in \dom p^{n+1}$
 $p^{n+1}_{\ordered{\gn_0, \dotsc, \gn_{l-1}}} \in i_{n, n+1}(D(A_n))$.
Hence
\begin{align*}
 i_{n+1,n+1+l-1}(p^{n+1})_{\ordered{\mc(p^{n+1}), \dotsc, 
        i_{n+1,n+l-1}(\mc(p^{n+1}))}} 
        \in i_{n,n+l}(D(A_n)).
\end{align*}
So $i_{n+1,\gw}(p^{n+1})_{\ordered{\mc(p^{n+1}), \dotsc, 
        i_{n+1,n+l-1}(\mc(p^{n+1}))}} 
        \in D(A)$.
Since $G$ is upward closed we get $G \intersect A \not= \emptyset$.
\end{proof}
\begin{theorem}
Assume GCH, $j \func V \to M \supset M^\gk$, $\crit(j) = \gk$,
	$\gen{j} \subset j(\gk)$.
Let $E$ be the extender derived from $j$.
Then there is a Cohen extension $V[G^*]$ in which there is 
$G$ which is $j_\gw(\PE)$-generic over
$M_\gw$.
\end{theorem}
\begin{proof}
Let $G^*$ be $\ordered{\PES, \leq^*}$-generic.
In $V[G^*]$ we set for each $n < \gw$
\begin{align*}
&  G_n = \setof {
		q \in j_n(\PE)
		} 
		{
		 f^q = f^{
			  j_{n}(p)_{
				\ordered{
					\mc(p), \dotsc, j_{n-1}(\mc(p))
				        }
				   }
			  },\ 
			p \in \PE,\ f^p \in G^*
		},
\\
& G = \setof {
		q \geq j_{n,\gw}(p)
              } 
		{n < \gw,\ 
	    		p \in G_n}.
\end{align*}		      
We show that $G$ is $j_\gw(\PE)$-generic over $M_\gw$.
Let $D_\gw \in M_\gw$ be a dense open subset of $j_\gw(\PE)$.
Let $D^*_\gw = \setof {f^p \in j_\gw(\PES)} {
  p \in j_\gw(\PE),\ 
 \exists l < \gw\ 
  \forall \ordered{\gn_0, \dotsc, \gn_{l-1}}\in \dom p\ 
       p_{\ordered{\gn_0, \dotsc, \gn_{l-1}}} \in D_\gw
}$. We take $g,h,\ga \in V$ such that 
\begin{align*}
& j_\gw(g)(\ga,\dotsc,j_{n-1}(\ga)) = D_\gw,
\\
& j_\gw(h)(\ga,\dotsc,j_{n-1}(\ga)) = D^*_\gw.
\end{align*}
Since $D^*_\gw$ is a dense open subset of $j_\gw(\PES)$ we have
\begin{align*}
X=\setof {\ordered{\gm_0, \dotsc, \gm_{n-1}}} {h(\gm_0, \dotsc, \gm_{n-1})
		\text{ is dense open in }\PES} \in E(\ga).
\end{align*}
Setting $D^* = \bigintersect_{\ordered{\gm_0, \dotsc, \gm_{n-1}} \in X}
		h(\gm_0, \dotsc, \gm_{n-1})$, we get that
$D^* \in V$ is a dense open subset of $\PES$. 
Using genericity we take
$f \in D^* \intersect G^*$.
By the construction of $D^*$ we have for each
$\ordered{\gm_0, \dotsc, \gm_{n-1}} \in X$ a condition
$q(\gm_0, \dotsc, \gm_{n-1}) \in h(\gm_0, \dotsc, \gm_{n-1})$ 
such that $f^{q(\gm_0, \dotsc, \gm_{n-1})} = f$.
This means that there are $X' \subseteq X$, $X' \in E(\ga)$, $l < \gw$
such that
\begin{multline*}
\forall \ordered{\gm_0, \dotsc, \gm_{n-1}} \in X'\ 
\forall \ordered{\gn_0, \dotsc, \gn_{l-1}} \in \dom q(\gm_0, \dotsc, \gm_{n-1})\ 
\\
q(\gm_0, \dotsc, \gm_{n-1})_{\ordered{\gn_0, \dotsc, \gn_{l-1}}} \in
	 g(\gm_0, \dotsc, \gm_{n-1}).
\end{multline*}
Set $q^* = j_n(q)(\ga, \dotsc, j_{n-1}(\ga))$.
Then $j_{n,\gw}(q^*)_{\ordered{\mc(q^*), \dotsc, 
				j_{n,n+l-1}(\mc(q^*))}} \in D_\gw$.
\end{proof}
\section{$\PE$-Forcing: The general case} \label{PEForcingGeneral}
We give here only the definition of the Prikry on Extender forcing 
for the general
case. We do not give proofs since they are a tedious
repetition of the proofs appearing in section \ref{PEForcing}.

Our assumptions are GCH and the existence of an elementary embedding
$j\func V \to M \supset M^\gl$, $\crit(j) = \gk$,
where $\gl$ is the minimal cardinal satisfying
$j(\gl) \supset \gen{j}$.
Let $E$ be the
extender derived from $j$. Recall $\dom E = j(\gl) \setminus j''\gl$.
\begin{definition}
Assume $T \subseteq [\gl]^{\upto \gw}$. 
For
        $\ordered{\gm_0, \dotsc, \gm_k},
                \ordered{\gn_0, \dotsc, \gn_n} \in T$
we define
        $\ordered{\gm_0, \dotsc, \gm_k} \ltT 
                \ordered{\gn_0, \dotsc, \gn_n}$
if
\begin{enumerate}
\item
        $k < n$.
\item
        $\ordered{\gm_0, \dotsc, \gm_k} =
                \ordered{\gn_0, \dotsc, \gn_k}$.
\end{enumerate}
Clearly $\ordered{T, \ltT}$ is a tree.
We always assume that $\ordered{} \in T$.
\end{definition}
\begin{definition}
Assume $T \subseteq [\gl]^{\upto \gw}$ is a tree. Then
\begin{enumerate}
\item
        $\Suc_T(\ordered{\gn_0, \dotsc, \gn_k}) =
                 \setof {\gn < \gl}
                {\ordered{\gn_0, \dotsc, \gn_k, \gn} \in T}$.
\item
        $\forall k<\gw$  $\Lev_k(T) = T \intersect [\gl]^k$.
\end{enumerate}
Note that $\ordered{} \in T$ implies $\:\Lev_0(T) = \set{\ordered{}}$.
\end{definition}
\begin{definition}
Assume $T \subseteq [\gl]^{\upto \gw}$ is a tree, 
        $\ordered{\gn_0,\dotsc, \gn_k} \in T$. Then
\begin{align*}
T_{\ordered{\gn_0, \dotsc, \gn_k}} = 
        \setof {\ordered{\gn_{k+1}, \dotsc, \gn_n} \in [\gl]^{\upto \gw} }
                {\ordered{\gn_0, \dotsc, \gn_k, \gn_{k+1}, \dotsc, \gn_n} \in T}.
\end{align*}
\end{definition}
\begin{definition}
Assume $T \subseteq [\gl]^{\upto \gw}$ is a tree, $A \in [\gl]^k$. Then
\begin{align*}
T\restricted A =
        \setof {\ordered{\gn_{1}, \dotsc, \gn_n} \in T}
                {n < \gw,\ \ordered{\gn_0, \dotsc, \gn_k} \in A}.
\end{align*}
\end{definition}
\begin{definition}
Assume $T^\gx \subseteq [\gl]^{\upto \gw}$ is a tree for all $\gx < \gl$.
Then $T = \bigintersect_{\gx < \gl} T^\gx$ is defined by induction on $k$ as:
\begin{enumerate}
\item
        $\Lev_0(T) = \set{\ordered{}}$.
\item
        $\ordered{\gn_0, \dotsc, \gn_k} \in T \implies$
                $\Suc_T({\ordered{\gn_0, \dotsc, \gn_k}}) =
                        \bigintersect_{\gx < \gl} 
                               \Suc_{T^\gx}(\ordered{\gn_0, \dotsc, \gn_k})$.
\end{enumerate}
\end{definition}
\begin{definition}
Let $F$ be a function such that $\dom F \subseteq [\gl]^{\upto \gw}$
is a tree and
        $\ordered{\gn_0,\dotsc, \gn_k} \in \dom F$. Then
        $F_{\ordered{\gn_0, \dotsc, \gn_k}}$ is a function such that
\begin{enumerate}
\item
        $\dom (F_{\ordered{\gn_0, \dotsc, \gn_k}}) = 
                (\dom F)_{\ordered{\gn_0, \dotsc, \gn_k}}$.
\item
        $F_{\ordered{\gn_0, \dotsc, \gn_k}} (\gn_{k+1}, \dotsc, \gn_n) = 
                F(\gn_0, \dotsc, \gn_k, \gn_{k+1}, \dotsc, \gn_n)$.
\end{enumerate}
\end{definition}
\begin{definition}
$T \subseteq [\gl]^{\upto \gw}$ is $E(\ga)$-tree if
        $\forall \ordered{\gn_0, \dotsc, \gn_k} \in T$
                $\Suc_T(\ordered{\gn_0, \dotsc, \gn_k})
                                \in E(\ga)$.
\end{definition}
We recall the definition of filter product in order to define powers
of $E(\ga)$.
\begin{definition}
We define powers of $E(\ga)$ by induction as follows:
\begin{enumerate}
\item
        For $k = 0$: $E^0(\ga) = \set{\emptyset}$.
\item
        For $k > 0$:
                $\forall A \subseteq [\gl]^k$ $A \in E^{k}(\ga) \iff$
        \begin{multline*}
        \setof{\ordered{\gn_0,\dotsc, \gn_{k-1}} \in [\gl]^{k-1}} 
        {
\\
        \setof {\gn_{k} \in \gl} 
                {\ordered{\gn_0, \dotsc, \gn_{k-1}, \gn_k} \in A} \in E(\ga)
        } \in E^{k-1}(\ga).
\end{multline*}
\end{enumerate}
Note that $E^1(\ga) = E(\ga)$. Recall that $E^k(\ga)$ is $\gk$-closed
ultrafilter on $[\gl]^k$.
\end{definition}
The following is straightforward.
\begin{proposition}
Assume $T \subseteq [\gl]^{\upto \gw},\  T^\gx \subseteq [\gl]^{\upto \gw}$ are
         $E(\ga)$-trees for all $\gx < \gm$, $\gm < \gk$. Then
\begin{enumerate}
\item
        $\forall k<\gw\ \Lev_k(T) \in E^k(\ga)$.
\item
        $\ordered{\gn_0, \dotsc, \gn_k} \in T \implies$
        $T_{\ordered{\gn_0, \dotsc, \gn_k}}$ is $E(\ga)$-tree.
\item
        $A \in E^k(\ga) \implies$
                $T \restricted A$ is $E(\ga)$-tree.
\item
        $\bigintersect_{\gx < \gm} T^{\gx}$ is $E(\ga)$-tree.
\end{enumerate}
\end{proposition}
We are ready to present the forcing notion:
\begin{definition}
A condition $p$ in $\PE$ is of the form
        $\ordered{f, \ga, F}$
where
\begin{enumerate}
\item
    $f\func d \to [\gl]^{\upto \gw}$ is such that
\begin{enumerate}
        \item
                $d \in [\dom E]^{\uptoeq \gl}$.
        \item
                $\forall \gb \in d$  $\ga \geqE \gb$.
\end{enumerate}
\item
        $F \func T \to [d]^{\upto \gl}$ is such that
\begin{enumerate}
        \item
                $T$ is $E(\ga)$-tree.
         \item 
                $\forall \ordered{\gn_0, \dotsc, \gn_k, \gn} \in T$
                        $F_{\ordered{\gn_0, \dotsc, \gn_k}}(\gn) \subseteq
                                        d$. 
        \item
                $\forall \ordered{\gn_0, \dotsc, \gn_k} \in T$
                        $j(F_{\ordered{\gn_0, \dotsc, \gn_k}})(\ga) = 
                                                                   j''d$.
        \item   
                $\forall \ordered{\gn_0, \dotsc, \gn_k, \gn} \in T$
                        $\gk \in F_{\ordered{\gn_0, \dotsc, \gn_k}}(\gn)$.
                Note: This condition implies $\gk \in d$!
        \item  
                $\forall \gb \in d$ $\forall \ordered{\gn_0, \dotsc, \gn_k} 
                                                        \in T\ $
                \begin{align*}
                f(\gb) \append \ordof {\gp_{\ga,\gb}(\gn_i)} {i \leq k,\,
                        \gb \in F(\gn_0, \dotsc, \gn_i)}
                        \in [\gl]^{\upto \gw}.
                \end{align*}
\end{enumerate}
\end{enumerate}
We write $\supp p$, $\mc(p)$, $f^p$, $F^p$, $\dom p$ for
         $d$,       $\ga$,    $f$,   $F$, $T$.
\end{definition}
%\fi
%
%
%
\begin{definition} 
Let $p, q \in \PE$. We say that $p$ is a Prikry extension of $q$
($p \leq^* q$ or $p \leq^0 q$) if
\begin{enumerate}
\item
        $\supp p \supseteq \supp q$.
\item
        $f^p \restricted \supp q = f^q$.
\item 
        $\dom p \subseteq \gp_{\mc(p), \mc(q)}^{-1} \dom q$.
\item
        
        $\forall \ordered{\gn_0,\dotsc,\gn_k, \gn} \in \dom p$
        $\forall \gb \in F^q \circ \gp_{\mc(p), \mc(q)}
                        {}_{\ordered{\gn_0, \dotsc, \gn_k}}(\gn)$
        \begin{align*}
		\gp_{\mc(p), \gb}(\gn) =
			\gp_{\mc(q), \gb}(
				\gp_{\mc(p), \mc(q)}(\gn)).
        \end{align*}
\item
        $\forall \ordered{\gn_0,\dotsc,\gn_k} \in \dom p$
                $F^p(\gn_0, \dotsc, \gn_k) \supseteq
                        F^q \circ \gp_{\mc(p),\mc(q)}
                         (\gn_0, \dotsc, \gn_k)$.
\item
        $\forall \ordered{\gn_0,\dotsc,\gn_k} \in \dom p$
                $F^p(\gn_0, \dotsc, \gn_k) \setminus
                        F^q \circ \gp_{\mc(p),\mc(q)}
                         (\gn_0, \dotsc, \gn_k) \subseteq 
                                \supp p \setminus \supp q$.
\end{enumerate}
\end{definition}
\begin{definition} 
Let $q \in \PE$ and $\ordered{\gn} \in \dom q$. 
We define 
$q_{\ordered{\gn}} \in \PE$ to be $p$ where
\begin{enumerate}
\item
        $\supp p = \supp q$.
\item
        $\forall \gb \in \supp p$
	$f^p(\gb)=$
        $\begin{cases}
        f^q(\gb) \append \ordered{\gp_{\mc(q), \gb}(\gn)}   &
                \gb \in F^q(\gn).
        \\
	f^q(\gb)   &
                \gb \notin F^q(\gn).
        \end{cases}$
\item
        $\mc(p) = \mc(q)$.
\item
        $F^p = F^q_{\ordered{\gn}}$.
\end{enumerate}
\end{definition}
When we write $q_{\ordered{\gn_0, \dotsc, \gn_k}}$ we mean
        $(\dotsb (q_{\ordered{\gn_0}})_{\ordered{\gn_1}}\dotsb)_
                        {\ordered{\gn_k}}$.
\begin{definition}
Let $p,q \in \PE$. We say that $p$ is a $1$-point extension of $q$
($p \leq^1 q$) if
there is $\ordered{\gn} \in \dom q$ such that
        $p \leq^* q_{\ordered{\gn}}$.
\end{definition}
%
%
%
% n-point extension
%
\begin{definition}
Let $p,q \in \PE$. We say that $p$ is an $n$-point extension of $q$
($p \leq^n q$) if there are $p^n, \dotsc, p^0$ such that
\begin{align*}
        p=p^n \leq^1 \dotsb \leq^1 p^0=q.
\end{align*}
\end{definition}
\begin{definition}
Let $p,q \in \PE$. We say that $p$ is an extension of $q$
($p \leq q$) if there is $n$ such that $p \leq^n q$.
\end{definition}
\begin{definition}
$\PES = \setof {f^p} {p \in \PE}$. The induced partial order $\leq^*$ 
on $\PES$ is $f \leq^* g \iff f \supseteq g$.
\end{definition}
\begin{definition}
Let $G$ be $\PE$-generic. Then
\begin{align*}
\forall \ga \in \dom E\ 
        G^\ga = \bigunion \setof {f^p(\ga)} {p \in G,\ \ga \in \supp p}.
\end{align*}
We write $\CN{G}^\ga$ for the $\PE$-name of $G^\ga$.
\end{definition}
\begin{proposition}
Let $G$ be $\PE$-generic. Then in $V[G]$:
\begin{enumerate}
\item $\ot G^\ga = \gw$.
\item $G^\ga$ is unbounded in $\gm$, where $\gm$ is the minimal ordinal
	such that 	$j(\gm) > \ga$.
\item $\ga \not= \gb \implies G^\ga \not= G^\gb$.
\end{enumerate}
\end{proposition}
\begin{claim}
Assume $p \in \PE$ and $\gs$ is a statement in the $\PE$-forcing
language. Then there is $p^* \leq^* p$ such that
	$p^* \decides_{\PE} \gs$.
\end{claim}
\begin{proposition}
$\ordered{\PE, \leq^*}$ is $\gk$-closed.
\end{proposition}
\begin{proposition}
Forcing with $\ordered{\PE, \leq^*}$ is the same as forcing with
$\ordered{\PES, \leq^*}$.
\end{proposition}
\begin{proposition}
$\PE$ satisfies $\gl^{++}$-cc.
\end{proposition}
\begin{claim}
$\PE$ is proper with regard to submodels of size $\gl$.
\end{claim}
\begin{theorem} \label{PrikryExtenderMainGeneral}
Assume $j \func V \to M \supset M^\gl$, $\crit(j) = \gk$,
$\gen{j} \subset j(\gl)$. Let $E$ be the extender derived from $j$,
and let $G$ be $\PE$-generic. Then in
$V[G]$:
\begin{enumerate}
\item
        $\forall \gm \in [\gk,\gl]_{\text{Reg}}$ $\cf \gm = \gw$.
\item
        $2^\gk = \power{j(\gl)}$.
\item
        No new bounded subsets are added to $\gk$.
\item
        All   cardinals outside of $(\gk,\gl]$ are preserved.
\end{enumerate}
\end{theorem}
%
%
%
%
%
%
%
%
%
%
%
%\bibliographystyle{plain}
%\bibliography{carmi}

\end{document}